\documentclass[a4paper,12pt]{compositio}

 \usepackage{amsfonts}
\usepackage{amscd}
\usepackage{amsmath, amssymb,latexsym}
\usepackage{amscd}
\usepackage{amssymb}
\usepackage{times}

 \usepackage[all] {xy}

\newtheorem{theo}{Theorem}[section]
\newtheorem{ex}[theo]{Example}
\newtheorem{prop}[theo]{Proposition}
\newtheorem{lem}[theo]{Lemma}
\newtheorem{cor}[theo]{Corollary}
\newtheorem{conj}[theo]{Conjecture}
\newtheorem{defi}[theo]{Definition}

\newtheorem{rema}[theo]{Remark}

\newtheorem*{maintheo}{{Theorem.}}

 \allowdisplaybreaks[4] 
\newcommand{\qeed}{\hfill\textrm{QED}\break\null}
\newenvironment{demo}{\noindent\textit{Proof.}~}{\qeed}

\def\bb{\mathbb}
       
\begin{document}
\title[Geometry of SSRs]{The geometry of special symplectic representations}

 \title[Special Symplectic Representations]{The geometry of special symplectic representations}

 \author{Marcus J. Slupinski }
 \email{marcus.slupinski@math.unistra.fr}
 \address{IRMA, Universit\'e de Strasbourg
 7 rue Ren\'e Descartes,
F-67084 Strasbourg Cedex France}
\author{Robert J. Stanton}
\email{stanton@math.ohio-state.edu}
 \address{Department of Mathematics, Ohio State 
University, 231 West 18th Avenue, 
Columbus OH 43210-1174}
\shortauthors{Slupinski-Stanton}



\begin{abstract} We show  there is a class of symplectic Lie algebra representations over any field of characteristic not $2$ or $3$  that have many    of the exceptional   algebraic and geometric properties  of both symmetric three forms in two dimensions  and   alternating three forms in six dimensions.  All nonzero orbits are coisotropic and the covariants satisfy relations generalising classical identities of  Eisenstein and Mathews. The main algebraic result is  that suitably generic elements of these representation spaces can be uniquely written  as the sum of two elements of a naturally defined Lagrangian subvariety.  We give universal explicit formulae for the summands and show how they lead to the existence of  geometric structure on appropriate subsets of  the representation space.
Over the real numbers this structure reduces to either a conic, special pseudo-K\" ahler metric or a conic, special para-K\" ahler metric.
\end{abstract}

\maketitle

\section{Introduction}
It has been known since the mid 19th century that symmetric three forms in two dimensions (binary cubics) possess  remarkable algebraic   properties. More recently \cite{Nigel} it was shown that real alternating three forms in six dimensions   also  have special algebraic properties. A common feature of these two spaces is that they are quite naturally symplectic vector spaces and  there is a natural choice of Lie algebra acting symplectically on them.  In the case of real three forms, Hitchin \cite{Nigel}  exploited this observation  extensively,  and although historically 
the symplectic aspect with regard to binary cubics has been largely  ignored, we  \cite{Sl-St1} showed that many of their important  properties can  be expressed in purely symplectic  terms.  The main  purpose of this paper is  to show that  there is a class of symplectic Lie algebra representations over any field $k$ of characteristic not $2$ or $3$  that have many    of the remarkable     properties  of both binary cubics  and   alternating three forms in six dimensions. These representations,
which we call  special symplectic representations  (SSR), are necessarily rare but   examples  include notably: 
\begin{itemize}
\item
$\mathfrak{sl}(2,k)$ acting on homogeneous polynomials of degree three in two variables ;
\item
$\mathfrak{sp}(6,k)$ acting on primitive alternating three forms in six dimensions; 
\item
$\mathfrak{sl}(6,k)$ acting on alternating three forms in six dimensions; 
\item
a $k$-form of $\mathfrak{so}(12,\bar{k})$ acting in a half-spinor representation  defined over $k$. 
\end{itemize}

\noindent The real vector space of alternating three forms in six dimensions is doubly interesting  because its special algebraic properties have  geometric implications. For example, certain of its open subsets are naturally  endowed with  \lq\lq special\rq\rq  differential  geometric structure (e.g. a special pseudo-K\" ahler metric  \cite{Nigel}). In the last section of the paper we show that the special algebraic properties of  all SSRs lead to special geometry so that, suitably interpreted, this structure exists for  special symplectic representations over a field $k$ of characteristic not $2$ or $3$.

 A special symplectic representation  is a symplectic Lie algebra representation  with extra structure. If $\mathfrak m$ denotes the Lie algebra and $V$ the representation space, this extra structure is an equivariant quadratic map $\mu:V\rightarrow \mathfrak m$  satisfying a constraint (cf Definition \ref{ssrdefi}). 
From $\mu$ and the symplectic form $\omega$ one can form two other symplectic covariants  $\Psi:V\rightarrow V$ and $Q:V\rightarrow k$, and these are the main technical tools of the paper. The three symplectic covariants generalise to any SSR the classical covariants of a binary cubic  defined by Eisenstein \cite{E}. 

Our first results describe orbit properties of SSRs or, more  precisely,  properties of  the vector space $\mathfrak m\cdot v$,  the tangent space to a group orbit through  $v\in V$ if  the action of $\mathfrak m$ is integrable. An unusual and important property is that $\mathfrak m\cdot v$  is coisotropic  if $v\not=0$,  and  we think this property may characterise SSRs. For binary cubics it is more or less evident  but seems not to have been known for other SSRs. 
Of particular interest  are generic orbits and minimal orbits.
We  show that  $\mathfrak m\cdot v$ is of codimension one if $Q(v)\not=0$ and that  $\mathfrak m\cdot v$ is Lagrangian iff $\mu(v)=0$ and $v\not=0$. In particular, an SSR is  a prehomogeneous vector space for the Lie algebra $k\times\mathfrak m$ and the open orbits provide examples of  noncommutative, completely integrable systems if $k=\bb R$ or $k=\bb C$. Although prehomogeneous vector spaces  for algebraic groups  have been  very widely studied  in the literature, there has not been a systematic investigation of their  local/global geometric structure.

The central  result of the paper is the Lagrangian decomposition theorem which involves all three symplectic covariants in an essential way. The existence of the following decomposition  for complex binary cubics was known to  Dickson \cite{Di}, and Hitchin \cite{Nigel} proved both  its existence and uniqueness   for suitably generic real alternating three forms in six dimensions.
\begin{maintheo} {\it Let $({\mathfrak m},V,\omega,B_{\mu})$ be an SSR and let $A\in V$. 

\noindent(a)  $Q(A)$ is a nonzero square in $k$ iff there exist $B,C\in V$ such that
$$
A=B+C,\quad \mu(B)=\mu(C)=0,\quad\omega(B,C)\not= 0.
$$

\noindent(b) $B$ and $C$ of (a) are unique up to permutation. In fact there is a square root $q$
 of $Q(A)$ such that  
 $$
 B=\frac{1}{2}(A+\frac{1}{q}\Psi(A)),\quad  C=\frac{1}{2}(A-\frac{1}{q}\Psi(A)),\quad q=- 3\omega(B,C).
 $$ 
 }
\end{maintheo}
\noindent In the case of the SSRs  given above   the zero set of $\mu$  is   composed of: 
\begin{itemize}
\item
binary cubics with  a triple root;
\item
decomposable primitive alternating three forms in six dimensions;
\item
decomposable alternating three forms in six dimensions; 
\item
pure spinors in twelve dimensions. 
\end{itemize}
For a general SSR, one may have to quadratically extend before there are  elements $A$ such that $Q(A)$ is a nonzero square or before the zero set of $\mu$ is nontrivial.

The    decomposition theorem has two important algebraic consequences. On the one hand it enables us to give a completely explicit description of the fibres of $\mu$ through generic points of $V$, and, on the other,  to prove that the  symplectic covariants of an SSR  satisfy a relation  which generalises the Eisenstein relation for the covariants of a binary cubic \cite{E}.

 Similarly, it has an important geometric consequence as it is also ultimately the source of   \lq\lq special" geometric structures on appropriate subsets of $V$.
This is explained in the last section of the paper where we are able to  formulate and  prove a generalisation of the special geometric result in \cite {Nigel} to the setting of arbitrary SSRs over fields of characteristic not 2 or 3. The generality of the field of coefficients necessitates a different approach to the formulation of special geometry than appears in the literature.

The authors would like to acknowledge the financial support of their respective institutions, notably the Universit\'e  de Strasbourg and Ohio State Univ. through its Math Research Institute, for extended visits that made possible this paper. In addition, RJS wants to acknowledge the support of the Max Planck Institut, Bonn, for an extended stay during which some of this research was done.

 \vskip0,2cm
{\it Throughout this paper we work over a field $k$ of characteristic not $2$ or $3$.}
\vskip0,2cm

\section{Special symplectic representations }
\subsection{ Definition and Background }
Let $\mathfrak m$ be a Lie algebra over $k$. The interesting representations of $\mathfrak m$ are usually those that possess a non-degenerate bilinear form compatible with the $\mathfrak m$ action. Much work has been done when the form is symmetric or Hermitian; however, the focus of this paper is on those that carry a special class of alternating form. The  definition of a special symplectic representation (SSR) we give below  is similar to the definition of 'special symplectic subalgebra'   in \cite{CS}  (cf Definition 2.3) but the difference is we are working over a more or less arbitrary field and do not assume the existence of an invariant nondegenerate symmetric bilinear form on $\mathfrak m$.
\begin{defi}\label{ssrdef}
 A  special symplectic representation $({\mathfrak m},V,\omega,B_{\mu})$  is a  faithful representation  of  ${\mathfrak m}$ on  a symplectic vector space $(V,\omega)$ together with an ${\mathfrak m}$-equivariant symmetric bilinear map 
$B_{\mu}:V\times V\rightarrow {\mathfrak m}$ such that for all $m\in{\mathfrak m}$ and $A,B,C\in V$,
 \begin{equation}\label{ssrdefi0}
\omega(m\cdot A,B)+\omega(A, m\cdot B)=0,
\end{equation}
\begin{equation}\label{ssrdefi}
2B_{\mu}(A,B)\cdot C-2B_{\mu}(A,C)\cdot B=
2\omega(B,C)A-\omega(A,B)C+\omega(A,C)B.
\end{equation}

\end{defi}

We write $\mu:V\rightarrow \mathfrak m$ for the quadratic map associated to  $B_{\mu}$  and set:
\begin{align}
\nonumber
&{\mathfrak m}_{\mu}={\rm Vect}<B_{\mu}(A,B)\in{\mathfrak m}:\, A,B\in V>,\\
\nonumber
&{\mathfrak m}^{\mu}=\{a\in \mathfrak {sp}(V,\omega):\, [a,B_{\mu}(v_1,v_2)] = 
B_{\mu}(a\cdot v_1,v_2)  + B_{\mu}(v_1,a\cdot v_2) \forall v_1,v_2\in V\} .
\end{align}
Clearly ${\mathfrak m}_{\mu}$ and ${\mathfrak m}^{\mu}$ are Lie subalgebras of $\mathfrak {sp}(V,\omega)$ and,
as $B_{\mu}$ is  ${\mathfrak m}$-equivariant,  ${\mathfrak m}_{\mu}$ is an ideal  in   ${\mathfrak m}^{\mu}$. Hence if ${\mathfrak m}'$ is any Lie algebra such that
$
{\mathfrak m}_{\mu}\lhd {\mathfrak m}' \subseteq {\mathfrak m}^{\mu},
$
then $({\mathfrak m}',V,\omega,B_{\mu})$ is also an SSR.

An SSR (like any symplectic representation)  has a  moment map $\tilde{\mu}:V\rightarrow{\mathfrak m}^*$  given by
\begin{equation}\label{classmom}
\tilde{\mu}(v)(m)=\omega(m\cdot v,v)\quad\forall v\in V,\forall m\in{\mathfrak m}.
\end{equation}
The existence of a non-degenerate, symmetric ${\mathfrak m}$-invariant bilnear form $(\phantom a,\phantom a)$ on ${\mathfrak m}$ such that
\begin{equation}\label{ssa}
(\mu(v),m)=\tilde{\mu}(v)(m)\quad\forall v\in V,\forall m\in{\mathfrak m},
\end{equation}
is equivalent to the image of ${\mathfrak m}$ in $\mathfrak {sp}(V,\omega)$  being a `special, symplectic subalgebra' in the sense of \cite{CS1}. In this case   it is easy to see that ${\mathfrak m}_{\mu}={\mathfrak m}$ and, if $k=\bb R$ or $k=\bb C$,  it can be shown that  ${\mathfrak m}_{\mu} ={\mathfrak m}^{\mu}$  (cf Proposition 2.7 in \cite{CS1}) .

The fundamental example of an SSR is the defining representation of a symplectic Lie algebra. If $(V,\omega)$ is a symplectic vector space, then $(\mathfrak {sp}(V,\omega),V,\omega,B_{\tau})$ is an SSR where 
\begin{equation}\label{inftransvec}
B_{\tau}(A,B)\cdot C=\frac{1}{2}\left(\omega(A,C)B+\omega(B,C)A\right),\quad\forall A,B,C\in V.
\end{equation}
With this observation it is informative to rewrite equation \eqref{ssrdefi} in the more natural equivalent  form:
\begin{equation}\label{newCS}
B_{\mu}(A,B)\cdot C  -B_{\mu}(A,C)\cdot B=B_{\tau}(A,B)(C)-B_{\tau}(A,C)(B).
\end{equation}

Some historical remarks are in order concerning Definition \ref{ssrdef}. To the best of the authors' knowledge, the main ingredient, equation \eqref{ssrdefi}, first appeared explicitly in the 1975 paper of W. Hein  (\cite{Hei}, Theorem 4, page 91). Following work of Faulkner, Freudenthal, Kantor, Koecher and Tits, he showed how one can construct a Lie algebra from a unital representation of a Jordan algebra   on a  triple system provided that these data are `admissible'  . Roughly speaking, equation  \eqref{ssrdefi} is the necessary and sufficient condition for admissibility of the data in the case of  a symplectic unital representation of a specific Jordan algebra. The results of Hein's paper are valid in characteristic not $2$ or $3$.

The same equation \eqref{ssrdefi} appears explicitly in the 2000 paper of L.J. Schwachh\"ofer (\cite{Sch1}, Theorem 6.1, page 304) but in a very different context. There,    it is a necessary, almost sufficient, condition for  a Lie algebra  ${\mathfrak m}$ over $\mathbb R$ to be the absolutely irreducible restricted holonomy Lie algebra of a {\it torsion free},  symplectic connection on some real manifold.

\subsection{ Symplectic covariants of an $SSR$}
Associated to a special symplectic representation we have the symplectic covariants. These are the  polynomial maps $\mu:V\rightarrow{\mathfrak m}$, $\Psi:V\rightarrow V$ and $Q:V\rightarrow k$ defined by
$$
\begin{cases}
\mu(A)=B_{\mu}(A,A),\\
\Psi(A)=\mu(A)\cdot A,\\
Q(A)=\frac{3}{2}\,\omega(A,\Psi(A))
\end{cases}
$$
and called, respectively, the quadratic, cubic and  quartic covariants of $({\mathfrak m},V,\omega,B_{\mu})$.  The  polar form of $\mu$ is $B_{\mu}$, and the polar forms of $\Psi$ and $Q$ are easily seen to be given by
\begin{align}\label{polcov}
\nonumber
B_{\Psi}(A,B,C)=&\frac{1}{3}\bigl(B_{\mu}(A,B)\cdot C+B_{\mu}(B,C)\cdot A+B_{\mu}(C,A)\cdot B\bigr),\\
\nonumber
B_{Q}(A,B,C,D)=&
\frac{3}{8}\bigl(\,\omega(A,B_{\Psi}(B,C,D))+\omega(B,B_{\Psi}(C,D,A))+\\
&\hskip1cm\omega(C,B_{\Psi}(D,A,B))+\omega(D,B_{\Psi}(A,B,C)\,\bigr).
\end{align}

The cubic and quartic covariants of a special symplectic representation are defined in terms of the quadratic covariant. A consequence of equation \eqref{newCS} is that any covariant  can be explicitly recovered from any other covariant.

\begin{prop}\label{covfromcov}
Let $({\mathfrak m},V,\omega,B_{\mu})$ be an SSR and let $A,B,C,D\in V$. Then:

\medskip
(i) $B_{\mu}(A,B)\cdot C-B_{\tau}(A,B)\cdot C=B_{\Psi}(A,B,C).$

\medskip
(ii) $\omega(A,B_{\mu}(B,C)\cdot D-B_{\tau}(B,C)\cdot D)=\frac{2}{3}B_{Q}(A,B,C,D).$

\medskip
(iii) $\omega(D,B_{\Psi}(A,B,C))=\frac{2}{3}B_{Q}(A,B,C,D).$
\end{prop}
\begin{demo} It is clear that  $B_{\mu}(A,B)\cdot C-B_{\tau}(A,B)\cdot C$ is invariant if we permute $A$ and $B$, and by  \eqref{newCS} it is also invariant if we permute $B$ and $C$. Hence 
$$
(A,B,C)\mapsto  B_{\mu}(A,B)\cdot C-B_{\tau}(A,B)\cdot C
$$
 defines a  symmetric trilinear map and, since $B_{\Psi}$ is also symmetric, to prove (i) it is sufficient to prove that for all $A\in V$,
 $$
 B_{\mu}(A,A)\cdot A-B_{\tau}(A,A)\cdot A=B_{\Psi}(A,A,A).
 $$
 The LHS is $ B_{\mu}(A,A)\cdot A=\Psi(A)$ which is  the RHS by \eqref{polcov}.
 
 Part (ii) is proved in a similar way. The LHS of (ii) is symmetric in $B,C,D$ by (i) and
 $$
  \omega(A,B_{\mu}(B,C)\cdot D-B_{\tau}(B,C)\cdot D)=
 \omega(D,B_{\mu}(B,C)\cdot A-B_{\tau}(B,C)\cdot A)
 $$
since $B_{\mu}(B,C)$ and  $B_{\tau}(B,C)$ are in $\mathfrak{sp}(V,\omega)$. Hence the LHS of (ii)  defines a symmetric quadrilinear form, so to prove (ii) it is sufficient to prove that for all $A\in V$,
 $$
 \omega(A,B_{\mu}(A,A)\cdot A-B_{\tau}(A,A)\cdot A)=\frac{2}{3}B_{Q}(A,A,A,A).
 $$
Since $B_{\tau}(A,A)\cdot A=0$, the LHS is $ \omega(A,B_{\mu}(A,A)\cdot A)=\omega(A,\Psi(A))=Q(A)$  which is  the RHS by \eqref{polcov}. Part (iii) follows directly from (i) and (ii).
\end{demo}
\begin{cor}\label{quadcovvan}
 Let $({\mathfrak m},V,\omega,B_{\mu})$ be an SSR.   If either the quartic covariant $Q$ or the cubic covariant $\Psi$ vanish identically, then $B_{\mu}=B_{\tau}$ and ${\mathfrak m}\cong \mathfrak {sp}(V,\omega)$.
 \end{cor}
\begin{demo} To simplify notation we identify ${\mathfrak m}$ with its image in $\mathfrak {sp}(V,\omega)$. If  the quartic covariant $Q$ (resp. the cubic covariant $\Psi$) vanishes identically
it follows from Proposition \ref{covfromcov}(ii)) (resp. Proposition \ref{covfromcov}(i))) that
$B_{\mu}=B_{\tau}$. Hence ${\mathfrak m}$ contains the Lie algebra generated by
all $B_{\tau}(A,B)$ (see \eqref{inftransvec}) and this is $\mathfrak {sp}(V,\omega)$.
\end{demo}

If $f:V\rightarrow k$ is a polynomial of degree not divisible by the characteristic of $k$, it is convenient to define  its `derivative' $df:V\rightarrow V^*$ by
$$
df_A(B)=({\rm deg}f)B_f(B,A,\cdots,A) \quad\forall A\in V,\forall B\in V,
$$
where $B_f$ is the  polar form of $f$.

\
\begin{cor} \label{smu}
 Let $({\mathfrak m},V,\omega,B_{\mu})$ be an SSR. For all $A\in V$, ${\rm Ker}\,d\mu_A=({\mathfrak m}_{\mu}\cdot A)^{\perp}$.
\end{cor}
\begin{demo} Interchanging $(A,B)$ and $(C,D)$  in the formula
$$
\omega(B_{\mu}(A,B)\cdot C,D)=-\frac{2}{3}B_Q(A,B,C,D)-\omega(B_{\tau}(A,B)\cdot C,D)
$$
(cf Proposition \ref{covfromcov}(ii)) and subtracting we get
\begin{align}
\nonumber
\omega(B_{\mu}(A,B)\cdot C,D)-&\omega(B_{\mu}(C,D)\cdot A,B)=\\
&\omega(B_{\tau}(A,B)\cdot C,D)-\omega(B_{\tau}(C,D)\cdot A,B).
\end{align}
Since $B_{\tau}(A,B)\cdot C=\frac{1}{2}\left(\omega(A,C)B+\omega(B,C)A\right)$, the RHS vanishes and hence
$$
\omega(B_{\mu}(A,B)\cdot C,D)=\omega(B_{\mu}(C,D)\cdot A,B)
$$
or, equivalently,
$$
\frac{1}{2}\omega(d\mu_A(B)\cdot C,D)=\omega(B_{\mu}(C,D)\cdot A,B).
$$
By definition, $(V,\omega)$ is a faithful symplectic representation of ${\mathfrak m}$ so it follows from this equation that
 $B\in{\rm Ker}\,d\mu_A$ iff $B\in ({\mathfrak m}_{\mu}\cdot A)^{\perp}$.
\end{demo}
\begin{rema} It is a well-known property of the moment map  \eqref{classmom} that for all $A\in V$, ${\rm Ker}\,d\tilde{\mu}_A=({\mathfrak m}\cdot A)^{\perp}$.
\end{rema}

\subsection{Examples}

The verification that a given symplectic representation of a Lie algebra satisfies the key equation (\ref{newCS}) to be an SSR  is often laborious and sometimes involves special dimension dependent identities. In this section we give examples of exceptional SSRs and   in the Appendix  examples of  SSRs which occur in infinite  families. The link with exceptional and
classical Lie algebras will be explained in the next section.

Many properties of  Example \ref{E6constr}  ($ \mathfrak {sl}(6,k)$ acting on $\Lambda^3({k^6}^*)$) were obtained by Hitchin in \cite{Nigel} without his observing explicitly the validity of (\ref{newCS}). 
From our point of view this SSR is obtained by `reducing' Example \ref{E7spinconstr} (a half-spin representation in dimension $12$) and Example \ref{F4constr} ($\mathfrak {sp}(3,k)$ acting on primitive $3$-forms in $\Lambda^3({k^6}^*)$) is obtained by `reducing' it.

\begin{ex} \label{E7spinconstr}
Let $(E,g)$ be a $12$-dimensional  vector space $E$ with a non-degenerate, hyperbolic, symmetric bilinear form  $g$. Let
 $\Sigma$ be either one of the $32$-dimensional irreducible half-spinor representations of $\mathfrak {so}(E,g)$. It was shown by E. Cartan \cite{Elie} that  there is a  (unique up to scaling) $\mathfrak {so}(E,g)$-invariant symplectic form $\omega$ on $\Sigma$ and, in terms of Clifford multiplication, he gave an explicit
 $\mathfrak {so}(E,g)$-equivariant surjection
$$
S^2(\Sigma)\rightarrow\Lambda^2(E)
$$
from the space of symmetric half-spinors to the space of bivectors. Composing this with 
an  $\mathfrak {so}(E,g)$-equivariant isomorphism $\Lambda^2(E)\cong \mathfrak {so}(E,g)$ we get an equivariant,  symmetric bilinear map  $B:\Sigma\times \Sigma\rightarrow \mathfrak {so}(E,g)$. After a delicate  spinor computation one can show that for fixed $\omega,$  a multiple of $B$, say $B_{\mu}$, satisfies \eqref{ssrdefi}.  Hence $(\Sigma,\mathfrak {so}(E,g),\omega,B_{\mu})$ is an SSR. Similarly, any other $k$-form of $\mathfrak {so}(E,g)\otimes_k\bar{k}$  whose half-spinors are defined over $k$ will  give rise to two SSRs.

For a hyperbolic quadratic form in dimension $4\, (mod\, 8)$, half-spinors are defined over $k$ and an invariant symplectic form $\omega$ and symmetric  map $B$  as above exist. However  it is only in dimensions $4$ and $12$ that  \eqref{ssrdefi} is  satisfied. In dimension $4$ the SSRs
we get  this way are both isomorphic to $\mathfrak{sl}(2,k)$ acting symplectically on $k^2$. 
\end{ex}

This SSR  gives rise two  other   exceptional SSRs  by  the following  `reduction' process.  Suppose that $({\mathfrak m},V,\omega,B_{\mu})$ is an SSR and ${\mathfrak h}\subseteq{\mathfrak m}$ is a Lie subalgebra such that $V^{{\mathfrak h}}$, the subspace annihilated by  $\mathfrak h$,  is a symplectic subspace. Let ${\mathfrak h}'$ be the commutant of ${\mathfrak h}$ in ${\mathfrak m}$.  Then if the action of ${\mathfrak h}'$   on 
$V^{{\mathfrak h}}$ is faithful,
$({\mathfrak h}',V^{{\mathfrak h}},\omega,B_{\mu})$ will also be an SSR since by equivariance, $B_{\mu}$ restricted to $V^{{\mathfrak h}}$ takes values in  ${\mathfrak h}'$ and already satisfies equation \eqref{ssrdefi}.

\begin{ex} \label{E6constr}
With the notation of  Example \ref{E7spinconstr}, we say that  an endomorphism $N:E\rightarrow E$ is a polarisation of $(E,g)$  if
$$
N^2=Id_{E},\quad N\in \mathfrak {so}(E,g).
$$
Then $E=E_1\oplus E_{-1}$ and the  eigenspaces $E_1,E_{-1}$ of $N$ are maximal isotropic. Furthermore, ${\rm ad}\,N:\mathfrak {so}(E,g)\rightarrow \mathfrak {so}(E,g)$ defines a 3-grading of $\mathfrak {so}(E,g)$, i.e.,
$$
\mathfrak {so}(E,g)={\mathfrak g}_{-2}\oplus{\mathfrak g}_0\oplus{\mathfrak g}_2
$$
where
$$
{\mathfrak g}_k=\{X\in \mathfrak {so}(E,g):\quad[N,X]=kX\}.
$$
The Lie algebra ${\mathfrak g}_0$ is isomorphic to $\mathfrak {gl}(6,k)$ since it preserves the decomposition $E=E_1\oplus E_{-1}$ and, as is easily seen,  restricting to $E_1$ (resp. $E_{-1}$) establishes an isomorphism ${\mathfrak g}_0\cong \mathfrak {gl}(E_1)$ (resp. ${\mathfrak g}_0\cong \mathfrak {gl}(E_{-1})).$    The derived algebra $\tilde{{\mathfrak g}_0}$ of ${\mathfrak g}_0$ is then  isomorphic to  $\mathfrak {sl}(6,k)$.

Polarisations of $(E,g)$ fall into two classes which can be distinguished by their spectrum in  the half-spinor representation $\Sigma$. In fact the action of  a polarisation $N$ on $\Sigma$ is diagonalisable and the set of eigenvalues is either $\{-3,-1,1,3\}$ or $\{-2,0,2\}$ so we have:
\begin{equation}\label{possspectra}
\Sigma=
\begin{cases}
\Sigma_{-3}\oplus\Sigma_{-1}\oplus\Sigma_{1}\oplus\Sigma_{3}
\quad\text{if spec}_{\Sigma}(N)=\{-3,-1,1,3\};\\
\Sigma_{-2}\oplus\Sigma_{0}\oplus\Sigma_{2}\quad\qquad\quad\text{if spec}_{\Sigma}(N)=\{-2,0,2\}.
\end{cases}
\end{equation}
The eigenspaces are stable under $\tilde{{\mathfrak g}_0}$ and in both cases  isomorphic to exterior powers of $E_1^*$ as  $\tilde{{\mathfrak g}_0}$-representations: 
\begin{equation}\label{partstates}
\Sigma_{k}\cong\Lambda^{3-k}(E_1^*)\quad \text{if}\quad k=-3,-2,\dots, 3.
\end{equation}
The symplectic form $\omega$ and  the map $B_{\mu}:\Sigma\times \Sigma\rightarrow \mathfrak {so}(E,g)$ of Example  \ref{E7spinconstr} are  $\mathfrak {so}(E,g)$-equivariant so for any polarisation:
$$
\omega(\Sigma_{m},\Sigma_{n})=0\,\text{if}\, m+n\not=0,\quad
{\mathfrak g}_{k}\cdot\Sigma_{m}\subseteq \Sigma_{k+m},\quad
B_{\mu}(\Sigma_{m},\Sigma_{n})\subseteq {\mathfrak g}_{m+n}.
$$
In particular, if ${\rm spec}_{\Sigma}(N)=\{-2,0,2\}$,  the restriction of $\omega$ to $\Sigma_{0}$,
the subspace annihilated by $N$, is nondegenerate. By reduction, it would follow that
$({\mathfrak g}_0,\Sigma_{0},\omega,B_{\mu})$ is an SSR if ${\mathfrak g}_0$ were to act  faithfully on $\Sigma_{0}$.  This is obviously false since $N\in{\mathfrak g}_0$  and $N$ acts trivially on $\Sigma_{0}$ by definition. However the derived algebra  $\tilde{{\mathfrak g}_0}$  acts faithfully by \eqref{partstates} and $B$ restricted to $\Sigma_{0}$ in fact takes its values in $\tilde{{\mathfrak g}_0}$ since there
 are no nontrivial symmetric scalar valued $\tilde{{\mathfrak g}_0}$-equivariant forms on $\Sigma_{0}$  ( a well-known property of  $sl(6,k)$acting on $\Lambda^3(k^6)$).  Hence $(\tilde{{\mathfrak g}_0},\Sigma_{0},\omega,B_{\mu})$ is an SSR. 
 
 To understand this SSR in terms intrinsic to the 6-dimensional vector space $E_1$,  we  identify $\Sigma_0$ with $\Lambda^{3}(E_1^*)$  and  the action of  $\tilde{{\mathfrak g}_0}$ on $\Sigma_0$ with  the natural action of ${\mathfrak sl}(E_1)$ on $\Lambda^3(E_1^*)$  (cf  \eqref{partstates} above). It can be shown that the symplectic form $\omega$  then gets identified with a symplectic form $\omega_1$ on $\Lambda^3(E_1^*)$ with the property
 $$
 \alpha\wedge \beta=\omega_1(\alpha,\beta){\bf vol}
 $$
 for some fixed  ${\bf vol}\in\Lambda^6(E_1^*)$, and the quadratic covariant $\mu$ with
 the map $\mu_1:\Lambda^3(E_1^*)\rightarrow {\mathfrak sl}(E_1)$ uniquely characterised by
 $$
\alpha\wedge i_e\alpha=i_{\mu_1(\alpha)(e)}{\bf vol}\qquad\forall e\in E_1,\,\forall\alpha\in \Lambda^3(E_1^*).
 $$
 
 Starting from these equations, Hitchin  \cite{Nigel}  examined this situation extensively .
  \end{ex}

By reduction of $(\tilde{{\mathfrak g}_0},\Sigma_{0},\omega,B)$ one  obtains  another exceptional SSR which can be described in terms of the `infinitesimal Hodge theory' of a 6-dimensional symplectic vector space.

\begin{ex} \label{F4constr}
With the notation of  Examples  \ref{E7spinconstr} and \ref{E6constr}, let $\Omega$ be a symplectic form  on the 6-dimensional vector space $E_1$ and let $\mathfrak {sp}(E_1,\Omega) $ be the corresponding symplectic Lie subalgebra of ${\mathfrak g}_0$. One can show  that $\mathfrak {sp}(E_1,\Omega)$ is its own double commutant  in $\mathfrak {so}(\Sigma,\omega)$ and that  its commutant $\mathfrak {h}$  in  $\mathfrak {so}(\Sigma,\omega)$ contains $N$ and is isomorphic to $\mathfrak {sl}(2,k)$. Decomposing $\Sigma_0$ under the action of $\mathfrak {sp}(E_1,\Omega)$  into irreducible components we get
$$
\Sigma_0={\bf 14}\oplus{\bf 6}
$$
where ${\bf 14}$ is the fixed point set of $\mathfrak {h}$ acting in $\Sigma$. By reduction
$(\mathfrak {sp}(E_1,\Omega),{\bf 14},\omega,B)$  is an SSR since $\mathfrak {sp}(E_1,\Omega)$ is simple and acts faithfully in any non-trivial  representation.

To understand this SSR in terms intrinsic to the 6-dimensional symplectic vector space $(E_1,\Omega)$, we proceed as in the previous example. Identifying $\Sigma$ with $\Lambda^{1}(E_1^*)\oplus\Lambda^{3}(E_1^*)\oplus\Lambda^{5}(E_1^*)$ (cf equations \eqref{possspectra} and \eqref{partstates} above), the action of  $\mathfrak {sp}(E_1,\Omega)$ on $\Sigma$ is identified with  the natural action of $\mathfrak {sp}(E_1,\Omega)$ on exterior forms of odd degree. One can show that  the action of $\mathfrak {h}$ is identified with the action of the `Hodge $\mathfrak {sl}(2,k)$' on  exterior forms of odd degree and hence the subspace  ${\bf 14}$  is identified with
the set of  $3$-forms which are annihilated by the Hodge operators, usually called primitive $3$-forms in Hodge theory and often denoted by $\Lambda^3_0(E_1^*)$.
\end{ex}
\begin{rema} 
The two pairs of Lie algebras $(\mathfrak {sp}(E_1),\Omega),\mathfrak {sl}(2,k)) $ and $({\mathfrak g}_0,kN)$  are  see-saw
dual pairs in $\mathfrak {so}(E,g)$, i.e.,  two pairs of mutual commutants in $\mathfrak {so}(E,g)$
with the property that $\mathfrak {sp}(E_1)\subseteq {\mathfrak g}_0$ and $\mathfrak {sl}(2,k)\supseteq kN$.
For either dual pair, the representation $\Sigma$ defines a Howe correspondence, i.e., the decomposition of $\Sigma$ into isotypic components for one member of the pair coincides with its decomposition into isotypic components for the other.
\end{rema}
The SSR of the next example was studied in great detail in \cite{Sl-St1}.
\begin{ex}\label{binarycubics} 
 Let $(E,\Omega)$ be a two-dimensional symplectic vector space. Since the characteristic of $k$
 is not $2$ or $3$,  we can identify $S^3(E^*)$, the set of  symmetric trilinear forms on $E$,  with  the set of cubic functions  on $E$.    This space has a  unique $sl(E)$-invariant symplectic structure $\omega$ such that
\begin{equation}\label{symp+ev0}
P(e)=\omega (P, {\tilde{e}}^3)\quad\forall e\in E,\,\forall P\in 
S^3(E^*),
\end{equation}
where \,${\tilde{}}:E\mapsto E^*$  is the unique isomorphism such that ${\tilde{e}}(f)=\Omega(e,f)$.
If  $\mu:V\rightarrow sl(E)$ is defined by
$$
{\rm Tr}_E(\mu(P)s)=\frac{1}{3}\omega(P,s\cdot P)\quad\forall s\in{\mathfrak m} ,
$$
the associated symmetric bilinear map $B_{\mu}$ satisfies   \eqref{ssrdefi} and so $(sl(E), S^3(E^*),\omega,B_{\mu})$  is an SSR.  Choosing a basis  $\{e_1,e_2\}$  of $E$ such that $\Omega(e_1,e_2)=1$ and  denoting by $\{x,y\}$ the dual basis of $E^*$, we can identify $S^3(E^*)$ with the set of homogeneous polynomials of degree three in $x$ and $y$. The 
above formulae then read:
\begin{align}
\nonumber
\omega(ax^3+3bx^2y+3cxy^2+dy^3,a'x^3+3b'x^2y&+3c'xy^2+d'y^3)\\
\nonumber
&=ad'-da'-3(bc'-cb').
\end{align}
\begin{equation}
\mu(ax^3+3bx^2y+3cxy^2+y^3)=
\begin{pmatrix}ad-bc&2(bd-c^2)\\2(b^2-ac)&-(ad-bc)
\end{pmatrix}
\end{equation}
and
\begin{equation}
Q(ax^3+3bx^2y+3cxy^2+y^3)=9\left(  (ad-bc)^2+4(bd-c^2)(b^2-ac)\right).
\end{equation}
Hence, up to  constants,   $Q(P)$   is  the discriminant of $P$ and $\mu(P)$ is the determinant of the Hessian of $P$ (modulo the $SL(2,k)$-equivariant identification of 
$\begin{pmatrix}
\alpha&\beta\\ \gamma&-\alpha
\end{pmatrix}\in \mathfrak {sl}(2,k)$ with the binary quadratic form $\beta x^2+2\alpha xy -\gamma y^2$).
\end{ex}

\subsection{Special symplectic representations - equivalent incarnations}
 In this section we associate a Lie algebra
${\mathfrak g}({\mathfrak m},V,\omega,B_{\mu})$ to an SSR $({\mathfrak m},V,\omega,B_{\mu})$ and characterise the simple Lie algebras obtained in this way by the existence of a certain type of grading.  Our approach is to show that an SSR is essentially the same thing as a vector space  with  a symplectic ternary product in the sense of Faulkner \cite{Fa}  and then use Faulkner's results.  In the context of real and complex special symplectic subalgebras,  similar results were proved by  Cahen-Schwachh\"ofer \cite{CS}  independently of \cite{Fa} but  using the existence of an invariant quadratic form on ${\mathfrak m}$ satisfying \eqref{ssa}. Although none of the results of this section will be used in the rest of the paper, we have included them to illustrate the r\^ole of SSRs in Lie theory.
\begin{defi}(\cite{Fa})
Let $V$ be a vector space. A Faulkner ternary product on V   is an antisymmetric bilinear form $ <\,,\,>:V\times V\rightarrow k$ and a trilinear map  $<\,,\,,\,>:V\times V\times V\rightarrow V$ such that:
\begin{align}
\nonumber
&(T1)\,<x,y,z>=<y,x,z>+<x,y>z\quad\forall x,y,z\in V;\\
\nonumber
&(T2)\,<x,y,z>=<x,z,y>+<y,z>x\quad\forall x,y,z\in V;\\
\nonumber
&(T3)\,<<x,y,z>,w>=<<x,y,w>,z>+<x,y><z,w>\\
\nonumber
&\hskip9,5cm\forall x,y,z,w\in V;\\
\nonumber
&(T4)\,<<x,y,z>,v,w>=<<x,v,w>,y,z>+<x,<y,v,w>,z>\\
\nonumber
&\hskip5cm+<x,y,<z,w,v>>\quad\forall x,y,z,v,w\in V.
\end{align}
The ternary product is called symplectic if $<\,,\,>$ is nondegenerate.
\end{defi}
The next proposition is the link between Faulkner ternary spaces and SSRs.
\begin{prop}  Let $V$ be a vector space and let   $ <\,,\,>:V\times V\rightarrow k$ be an antisymmetric bilinear form. Let $<\,,\,,\,>:V\times V\times V\rightarrow V$ and $B:V\times V\rightarrow {\rm End}(V)$ be respectively a trilinear and a symmetric bilinear map such that
\begin{equation} \label{CS=Fa}
B(x,y)\cdot z=\frac{1}{2}<x,y>z-<z,x,y>\quad\forall x,y,z\in V.
\end{equation} 
Then $(<\,,\,>,<\,,\,,\,>)$ is a Faulkner ternary product on $V$ iff $B$ satisfies
\begin{align}
\nonumber
&(B1)\quad B(x,y)\cdot z=B(x,z)\cdot y+<y,z>x-\frac{1}{2}<z,x>y+\frac{1}{2}<y,x>z\\
&\nonumber\hskip10cm\forall  x,y,z\in V;\\
\nonumber
&(B2)\quad B(x,y)=B(y,x)\quad\forall x,y\in V;\\
\nonumber
&(B3)\quad <B(x,y)\cdot u,v>+<u,B(x,y)\cdot v)=0\quad\forall x,y,u,v\in V;\\
\nonumber
&(B4)\quad [B(x,y),B(u,v)]=B(B(x,y)\cdot u,v)+B(u,B(x,y)\cdot v)\\
\nonumber
&\hskip9cm\quad\forall x,y,u,v,\in V.
\end{align}
\end{prop}
\begin{demo} It is immediate that (T1) is equivalent to (B1) and that (T2) is equivalent to (B2).
Substituting \eqref{CS=Fa} in (T3) gives
\begin{align}
\nonumber
\frac{1}{2}<y,z><x,w>-&<B(y,z)\cdot x,w>=\\
\nonumber
\frac{1}{2}<y,w><x,z>-&<B(y,w)\cdot x,z>+<x,y><z,w>
\end{align}
and then using (B1) to replace $B(y,z)\cdot x$ and $B(y,w)\cdot x$ by $B(y,x)\cdot z$ and $B(y,x)\cdot w$ respectively, this reduces to (B3). Hence (T1), (T2) and (T3) are equivalent to (B1), (B2) and (B3).

Substituting \eqref{CS=Fa} in the RHS of (T4) gives
\begin{align}
\nonumber
\frac{1}{2}<v,w><x,y,z>-&<B(v,w)\cdot x,y,z>\\
\nonumber
&-<x,B(v,w)\cdot y,z>-<x,y,B(v,w)\cdot z>.
\end{align}
and substituting \eqref{CS=Fa} in the LHS of (T4) gives
$$
\frac{1}{2}<v,w><x,y,z>-B(v,w)\cdot<x,y,z>.
$$
Hence (T4) implies
\begin{align}
\nonumber
B(v,w)\cdot<x,y,z>=&<B(v,w)\cdot x,y,z>\\
\nonumber
+&<x,B(v,w)\cdot y,z>+<x,y,B(v,w)\cdot z>,
\end{align}
i.e., (T4) implies that $<\,,\,,\,>$ is $B(v,w)$-equivariant for all $v,w\in V$ and the converse holds too by reversing the substitutions. From  \eqref{CS=Fa} it is clear that $<\,,\,,\,>$ is $B(v,w)$-equivariant for all $v,w\in V$ iff (B4).
\end{demo}

It follows from the proposition that if $(<\,,\,>,<\,,\,,\,>)$ is a symplectic Faulkner ternary product on $V$, then $({\mathfrak m},V,<\,,\,>,B)$ is an SSR where $\mathfrak {m}$ is  any `intermediate' Lie algebra. By `intermediate' we mean a Lie subalgebra of $\mathfrak {sp}(V,<\,,\,>)$ which
contains the Lie algebra generated by the image of $B$ and which is contained in the subalgebra  for which  $B$ is an equivariant map. Conversely,
if $({\mathfrak m},V,\omega,B_{\mu})$ is an SSR then $(\omega,<\, ,\,,\,>_{\mu})$ is a symplectic Faulkner ternary product on $V$ for which $\mathfrak {m}$ is intermediate  if we define $<\, ,\,,\,>_{\mu}$ from $B_{\mu}$ by equation \eqref{CS=Fa}.

We now recall the main results of \cite{Fa}.
Given a ternary product space  $(V,<\,,\,>,<\,,\,,\,>)$ and an  intermediate Lie algebra $\mathfrak {m}$,  Faulkner  \cite{Fa} defined a Lie bracket on the vector space
$$
\mathfrak {g}(\mathfrak {m}, V,<\,,\,>,<\,,\,,\,>)=\mathfrak {m}\oplus\mathfrak{sl}(2,k)\oplus V\otimes k^2
$$
extending the Lie bracket of $\mathfrak {m}\oplus\mathfrak{sl}(2,k)$ and the  action of $\mathfrak {m}\oplus\mathfrak{sl}(2,k)$ on $V\otimes k^2$.  
\begin{theo} \label{Fa1}
(\cite{Fa} Theorem 1) Let $(V,<\,,\,>,<\,,\,,\,>)$ be a ternary product space. Then $\mathfrak {g}(\mathfrak {m}, V,<\,,\,>,<\,,\,,\,>)$ is simple iff $<\,,\,>$ is nondegenerate and ${\mathfrak m}={\rm Im}\,B$.
\end{theo}

Faulkner also gave  a characterization of the Lie algebras obtained by his construction. In \cite {Ka} one can find a proof, under assumptions on $k$, that various  related constructions by Hein, Faulkner, Allison and Freudenthal  produce isomorphic Lie algebras.

 \begin{defi} \label{grading}
Let $\mathfrak{g}$ be a 
Lie algebra over $k$. A Heisenberg grading operator of $\mathfrak{g}$ is an element $H\in \mathfrak{g}$ such that

\medskip
(a) ad H is diagonalizable with eigenvalues $-2,-1,0,1,2$;

\medskip
(b) The eigenspaces corresponding to $\pm 2$ satisfy $\textrm{dim}(\mathfrak{g}_{2})$= dim$(\mathfrak{g}_{-2})=1$;

\medskip
(c) There exist $E\in\mathfrak{g}_{2}$ and $F\in\mathfrak{g}_{-2}$ with $\{E,H,F\}$ a standard 
$sl_2$-triple.

\medskip
\noindent We say $H$ is a Faulkner grading operator if $H$ also satisfies

\medskip
(d) The commutant ${\mathfrak m}$ of $\mathfrak{s} = <E,H,F>$ contains no nonzero ideals of $\mathfrak{g}$.
\end{defi}
Note that $-2,-1,0,1,2$ are distinct elements of $k$ since char$(k)\not=2,3$. 
For $n$ in $k$ we set 
$$\mathfrak{g_n} = \{Z\in \mathfrak{g}:[H,Z] = n Z\}.$$
Then,
$$\mathfrak{g}= \mathfrak{g}_{-2}\oplus 
\mathfrak{g}_{-1}\oplus\mathfrak{g}_{0}\oplus\mathfrak{g}_{1}\oplus\mathfrak{g}_{2},$$
and we call this the Heisenberg grading of $\mathfrak{g}$ associated to $H$.  In the literature such gradings are also known as `gradings of the second type', '$2$-gradings' or `contact gradings' but to our knowledge they first appeared in \cite{Fa}.  Heisenberg gradings and Faulkner gradings are the same thing if $\mathfrak{g}$ is simple.
\begin{theo} \label{Fa2}
(\cite{Fa} Theorem 2) 
Let $\mathfrak{g}$ be a Lie algebra over $k$. There exists a ternary product space $(V,<\,,\,>,<\,,\,,\,>)$ and intermediate Lie algebra $\mathfrak{m}$ such that $\mathfrak{g}\cong \mathfrak {g}(\mathfrak {m}, V,<\,,\,>,<\,,\,,\,>)$  iff there exists a Faulkner grading operator $H\in\mathfrak{g}$. 
\end{theo}

As we saw above, an SSR is equivalent to a symplectic ternary product space with intermediate Lie algebra.   Given an SSR $({\mathfrak m},V,\omega,B_{\mu})$ we define ${\mathfrak g}({\mathfrak m},V,\omega,B_{\mu})$ to be the Lie algebra $\mathfrak {g}({\mathfrak m},V,\omega,<\, ,\,,\,>_{\mu})$. 
\begin{theo}\label{SSRandsimpleLiealg}
(i) Let $({\mathfrak m},V,\omega,B_{\mu})$ be an SSR. Then ${\mathfrak g}({\mathfrak m},V,\omega,B_{\mu})$ is simple iff ${\mathfrak m}={\mathfrak m}_{\mu}.$

\noindent(ii) Let $\mathfrak{g}$ be a simple Lie algebra over $k$.  There exists an SSR $({\mathfrak m},V,\omega,B_{\mu})$ such that $\mathfrak{g}\cong {\mathfrak g}({\mathfrak m},V,\omega,B_{\mu})$  iff there exists a Heisenberg grading operator $H\in\mathfrak{g}$. 
\end{theo}
\begin{demo} This is immediate from Theorems \ref{Fa1} and  \ref{Fa2}.
\end{demo}
\begin{rema}   
Let $(\mathfrak{g},H,E,F)$ be a simple Heisenberg graded Lie algebra.  The authors show in \cite{Sl-St2} that $\mathfrak{g}\cong {\mathfrak g}({\mathfrak m},V,\omega,B_{\mu})$ where:
\begin{itemize}
\item $\mathfrak {m}$ is the commutant of $\mathfrak{s} = <E,H,F>$ in $\mathfrak{g}$;
\item $V=\mathfrak{g}_1$;
\item $\omega$ is defined by $[v_1,v_2]=\omega(v_1,v_2)E$ for all $v_1,v_2\in V$;
\item $B_{\mu}(v_1,v_2)=-\frac{1}{2}{\rm ad} (v_1)\circ{\rm ad} (v_2)(F)-\frac{1}{2}{\rm ad} (v_2)\circ{\rm ad} (v_1)(F)$ for all $v_1,v_2\in V$.
\end{itemize}
\end{rema}
Theorem \ref{SSRandsimpleLiealg} can be read in two ways: either as a means of constructing simple Lie algebras over $k$  (this was Faulkner's motivation) or as a means of constructing examples of SSRs. For the second point of view one needs to find examples of simple Lie algebras which have Heisenberg gradings.  Any finite-dimensional simple complex Lie algebra other than $A_1$ has 
 a unique (up to automorphism) Heisenberg grading. Indeed all such gradings can be obtained as follows:  take $E$ an element in the 
minimal nilpotent orbit in $\mathfrak{g}$ and $\{E,H,F\}$ a Jacobson-Morozov triple; by standard root theory the eigenvalues 
of ${\rm ad}H$ are $\{0,\pm 1,\pm 2\}$ and the eigenspaces of ${\rm ad}H$ define a simple Heisenberg grading. Cheng \cite{Ch} gave a complete classification
of finite-dimensional  simple real  Lie algebras admitting a Heisenberg grading. For general $k$ the authors are not aware of a classification. 

To end this section we give the list of Lie algebras associated to the exceptional SSRs of the previous section:
\begin{align}
\nonumber
\mathfrak {g}(\mathfrak {so}(V,g) ,\Sigma,\omega,B_{\mu})&\cong \mathfrak e_7\quad\text{(Example \ref{E7spinconstr}),}\\
\nonumber
\mathfrak {g}(\mathfrak {sl}(E_1) ,\Lambda^3(E_1^*),\omega,B_{\mu})&\cong \mathfrak e_6
\quad\text{(Example \ref{E6constr}),}\\
\nonumber
\mathfrak {g}(\mathfrak {sp}(E_1,\omega) ,\Lambda^3_0(E_1^*),\omega,B_{\mu})&\cong \mathfrak f_4\quad \text{(Example \ref{F4constr}),}\\
\nonumber
\mathfrak {g}(\mathfrak {sl}(E) ,S^3(E^*),\omega,B_{\mu})&\cong \mathfrak g_2\quad\text{(Example \ref{binarycubics}).}
\end{align}
The Lie algebras on the RHS are all split. Note that the two SSRs corresponding to the two non-isomorphic half-spin representations of $\mathfrak {so}(V,g)$  (Example \ref{E7spinconstr}) give rise to isomorphic Lie algebras. It is known (\cite{Fa} Example 2) that split $\mathfrak e_8$  is associated to the SSR corresponding to the irreducible $56$-dimensional representation of split $\mathfrak e_7$ with appropriate symplectic form and ternary product. However we did not include this example in the previous section as we know of no `geometric' description of this representation.

\section{Orbit properties of an $SSR$}

\subsection{Coisotropy}

An $SSR$  $({\mathfrak m},V,\omega,B_{\mu})$  is `special' in the algebraic sense that $B_{\mu}$ satisfies  algebraic constraints. In this section  we show that it is also `special' in the  geometric  sense that ${\mathfrak m}\cdot A$, the formal tangent space to the orbit through $A$, is coisotropic if $A\not=0$. If the action of ${\mathfrak m}$ on $V$ can be integrated to a group action the corresponding group orbits will then be  coisotropic but  the authors do not know whether this can be done over an arbitrary field of characteristic not 2 or 3.
\begin{theo} \label{orbcoiso} 
Let $({\mathfrak m},V,\omega,B_{\mu})$ be an SSR and let $A\in V\setminus\{0\}$. Then   ${\mathfrak m}\cdot A$ is coisotropic.
\end{theo}
\begin{demo} It is sufficient to show that ${\mathfrak m}_{\mu}\cdot A$ is coisotropic  since ${\mathfrak m}_{\mu}\cdot A\subseteq{\mathfrak m}\cdot A$.
Recall  that ${\mathfrak m}_{\mu}\cdot A$ is the symplectic orthogonal of ${\rm Ker}\,d\mu_A$ by Proposition \ref{smu}, so to prove ${\mathfrak m}_{\mu}\cdot A$ is coisotropic it suffices to show that ${\rm Ker}\,d\mu_A$ is isotropic, i.e., that $\omega(B,C)=0$ if $B,C\in{\rm Ker}\,d\mu_A$. 

Substituting
$$
d\mu_A(B)=2B_{\mu}(A,B)
$$
into  \eqref{ssrdefi}  gives
$$
d\mu_A(B)\cdot C-d\mu_A(C)\cdot B=2\omega(B,C)A-\omega(A,B)C+\omega(A,C)B
$$
for all $A,B,C\in V$. If $B,C\in{\rm Ker}\,d\mu_A$ this implies that
\begin{equation}\label{CSmom}
2\omega(B,C)A-\omega(A,B)C+\omega(A,C)B=0
\end{equation}
and hence, contracting with $B$ and $C$ successively, that
$$
3\omega(B,C)\omega(A,B)=3\omega(A,C)\omega(B,C)=0.
$$
From this it follows (since char$(k)\not=3$) that either $\omega(B,C)=0$ or $\omega(A,B)=\omega(A,C)=0$. But
by equation \eqref{CSmom}, $\omega(A,B)=\omega(A,C)=0$ implies (since char$(k)\not=2$) that $\omega(B,C)=0$ (since $A\not=0$)
and hence in all cases we have $\omega(B,C)=0$.  This proves the result.
 \end{demo}
 
 The property that all nontrivial orbits are coisotropic  is a very strong constraint on a symplectic representation. If $k=\bb C$,  finite-dimensional symplectic representations such that only the generic orbit is coisotropic   have been classified  (cf Knopp). It seems reasonable to conjecture
 
  \begin{conj}   Let $({\mathfrak m},V,\omega)$ be a symplectic representation of the Lie algebra ${\mathfrak m}$. Then ${\mathfrak m}\cdot A$ is coisotropic for all nonzero $A$ in $V$ iff there exists an ${\mathfrak m}$-equivariant symmetric bilinear map $B_{\mu}:V\times V\rightarrow {\mathfrak m}$ which  satisfies  \eqref{ssrdefi}.
 \end{conj}

\subsection{Generic orbits and minimal orbits} 
For certain $A\in V$ one can give more precise information on  ${\mathfrak m}\cdot A$
 using the following characterisation of the vanishing set of $Q$. (H. Rubenthaler informs us that he has a similar characterisation for certain prehomogeneous vector spaces associated to $3$-graded  Lie algebras.)
\begin{lem} \label{vanQchar}
Let $({\mathfrak m},V,\omega,B_{\mu})$ be an SSR and let
$A\in V$. Then
\medskip
$Q(A)=0$ iff  $A\in{\mathfrak m}\cdot A$.
\end{lem}
\begin{demo}
To prove this we show the sequence of implications:  $A\in{\mathfrak m}\cdot A$ $\Rightarrow$ $Q(A)=0$ $\Rightarrow$ $A\in{\mathfrak m}_{\mu}\cdot A$.
If $A\in{\mathfrak m}\cdot A$ there exists  $s\in{\mathfrak m}$ such that
$A=s\cdot A$. By Euler's formula and ${\mathfrak m}$-equivariance of $Q$, we get
$$
Q(A)=\frac{1}{4}\,dQ_A(A)=\frac{1}{4}\,dQ_A(s\cdot A)=0.
$$
Hence $A\in{\mathfrak m}\cdot A$ implies $Q(A)=0.$

To prove the second implication, suppose that $Q(A)=\frac{3}{2}\omega(A,\Psi(A))=0$. To show that $A\in{\mathfrak m}_{\mu}\cdot A$, it is sufficient to show that $\omega(A,B)=0$ for all $B\in{\rm Ker}{d\mu_A}$
since ${\mathfrak m}_{\mu}\cdot A={\rm Ker}{d\mu_A}^{\perp}$. 

Let $B\in{\rm Ker}\,d\mu_A$. Then for any $C\in V$, equation \eqref{ssrdefi} implies
\begin{align}
\nonumber
-2B_{\mu}(A,C)\cdot B&=
2\omega(B,C)A-\omega(A,B)C+\omega(A,C)B,\\
\nonumber
-2B_{\mu}(B,C)\cdot A&=
2\omega(A,C)B-\omega(B,A)C+\omega(B,C)C.
\end{align}
Taking $C=\Psi(A)$ and  $C=A$  respectively in the first equation, and $C=B$ in the second gives
the system
\begin{align}
\label{s1}
0&=2\omega(B,\Psi(A))A-\omega(A,B)\Psi(A),\\
\label{s2}
-2\mu(A)\cdot B&=-3\omega(A,B)A,\\
\label{s3}
-2\mu(B)\cdot A&=-3\omega(B,A)B,
\end{align}
where  to get the LHS of \eqref{s1} we have  used 
$$
B_{\mu}(A,\Psi(A))=B_{\mu}(A,\mu(A)\cdot A)=\frac{1}{2}[\mu(A),\mu(A)]=0.
$$
It follows from \eqref{s1} that
$$
\omega(A,B)\Psi(A)=2\omega(B,\Psi(A))A=2\omega(B,\mu(A)\cdot A)A
=-2\omega(\mu(A)\cdot B, A)A
$$
and substituting \eqref{s2} in this gives
$$
\omega(A,B)\Psi(A)=-3\omega(\omega(A,B)A,A)=0.
$$
Suppose for a contradiction that $\omega(A,B)\not=0$. Then $\Psi(A)=\mu(A)\cdot A=0$ and
$$
[\mu(A),\mu(B)]\cdot A=\mu(A)\cdot(\mu(B)\cdot A)=-\frac{9}{4}\omega(A,B)^2A
$$
using \eqref{s3} and then \eqref{s2}. Since  $\omega(A,B)\not=0$ this implies $A\in{\mathfrak m}_{\mu}\cdot A={\rm Ker}{d\mu_A}^{\perp}$ and therefore $\omega(A,B)=0$ which is a contradiction.  Thus $\omega(A,B)=0$ for all $B\in{\rm Ker}{d\mu_A}$ and hence $A\in{\mathfrak m}_{\mu}\cdot A$.
\end{demo}
 
 Generic and Lagrangian `orbits' are described by the following
\begin{prop} \label{skychara} 
Let $({\mathfrak m},V,\omega,B_{\mu})$ be an SSR and let $A\in V$.
\vskip0,1cm
\noindent(i)  $Q(A)\not=0$ iff $V={\mathfrak m}\cdot A\,\,\oplus<A>$.  In this case
${\rm Ker}\,{d\mu_A}=<\Psi(A)>$ and ${\rm Ker}\,{dQ_A}={\mathfrak m}\cdot A$. 
\vskip0,1cm
\noindent(ii) Suppose ${\rm dim}\,V> 2$ and  $A\not=0$. Then ${\mathfrak m}\cdot A$ is Lagrangian iff  $\mu(A)=0$.
\end{prop}
\begin{demo} (i): By Lemma \ref{vanQchar}, $V={\mathfrak m}\cdot A\,\,\oplus<A>$ implies
$Q(A)\not=0$. To prove the converse suppose $Q(A)\not=0$.
Since $Q(A)=\frac{3}{2}\omega(A,\Psi(A))$,  we have $\Psi(A)\not=0$ and by ${\mathfrak m}$-equivariance, 
$$
d\mu_A(\Psi(A))=2B_{\mu}(A,\Psi(A))=[\mu(A),\mu(A)]=0.
$$
Hence ${\rm Ker}{d\mu}_A$ is of dimension at least one. Let $B\in {\rm Ker}{d\mu}_A$. From equation {\eqref{ssrdefi} } with $C=\Psi(A)$, we get
$$
0=
2\omega(B,\Psi(A))A-\omega(A,B)\Psi(A)+\omega(A,\Psi(A))B
$$
and so $B$ is a linear combination of $A$ and $\Psi(A)$, say
$
B=aA+b\Psi(A).
$
But 
$$
0=d\mu_A(B)=ad\mu_A(A)+bd\mu_A(\Psi(A))=2a\mu(A)
$$
 and since $\mu(A)\not=0$,  this implies $a=0$ and $B\in<\Psi(A)>$. This proves that 
 ${\rm Ker}\,{d\mu_A}=<\Psi(A)>$ and it  follows that $V={\mathfrak m}\cdot A\,\,\oplus<A>$ since  this  sum   is direct (cf Lemma \ref{vanQchar}) and codim$({\mathfrak m}\cdot A)\le{\rm codim}({\mathfrak m}_{\mu}\cdot A)={\rm dim}\,{\rm Ker}\,{d\mu_A}=1$. The last part of (i) follows since 
  ${\mathfrak m}\cdot A\subseteq {\rm Ker}\,{dQ_A}$ (by ${\mathfrak m}$-invariance of  $Q$) and both spaces are of codimension one.

To prove (ii) suppose first that ${\mathfrak m}\cdot A$ is Lagrangian. Then  we must have
$Q(A)=0$ for if not, ${\mathfrak m}\cdot A$ is of codimension one by (i) and cannot be Lagrangian if ${\rm dim}{\,\mathfrak g_{1}}> 2$. By Proposition \ref{vanQchar}, 
this means $A\in{\mathfrak m}\cdot A$.
But ${\mathfrak m}\cdot A$ is Lagrangian so
$$
A\in {\mathfrak m}\cdot A=({\mathfrak m}\cdot A)^{\perp}\subseteq ({\mathfrak m}_{\mu}\cdot A)^{\perp}={\rm Ker}\,d\mu_A
$$ 
and hence $0=d\mu_A(A)=2\mu(A)$.

To prove  implication in the opposite direction, suppose $\mu(A)=0$. Then $\mathfrak{m}$-equivariance of $\mu$ implies
$$
d\mu_A(m\cdot A)=[m, \mu(A)]=0\quad\forall m\in\mathfrak{m}, 
$$
and so ${\mathfrak m}_{\mu}\cdot A\subseteq{\mathfrak m}\cdot A\subseteq{\rm Ker}\,d\mu_A=({\mathfrak m}_{\mu}\cdot A)^{\perp}$.  This means   ${\mathfrak m}_{\mu}\cdot A$ is isotropic but  by  Theorem  \ref{orbcoiso}, it is also
coisotropic if $A\not=0$. Hence ${\mathfrak m}_{\mu}\cdot A$ is Lagrangian and ${\mathfrak m}\cdot A={\mathfrak m}_{\mu}\cdot A$ is Lagrangian.
\end{demo}
\begin{rema} If ${\rm dim}{V}=2$  and $A\not=0$, one can show that either
$Q(A)\not=0$ or $\mu(A)=0$ and that in both cases ${\mathfrak m}\cdot A$ is of dimension one  and so Lagrangian.
\end{rema}
\begin{rema} 
A prehomogeneous vector space(PV)  is a representation of a group $G$ on a finite-dimensional vector space $V$  such that $G$ has a Zariski open  orbit in V.
By Proposition \ref{skychara}(i), taking the product of  the actions of $k$ and $\mathfrak{m}$ on $V$ gives an action of $\hat{\mathfrak{m}}=k\oplus\mathfrak{m}$ on $V$ such that $\hat{\mathfrak{m}}\cdot A=V$ if $Q(A)\not= 0$. Hence $(\hat{\mathfrak{m}},V)$ satisfies an  infinitesimal analogue of the PV condition.
\end{rema}
If  $k=\bb R$, we can exponentiate the action of ${\mathfrak m}$ on $V$ and  the  orbits of the corresponding group $M$  define a codimension one foliation ${\mathcal G}$ of the set of generic points ${\mathcal O}=\{A\in V: \,Q(A)\not=0\}$ by Proposition \ref{skychara}(i). The fibres of $\mu$  define a one-dimensional  foliation ${\mathcal F}$ of ${\mathcal O}$ also by Proposition \ref{skychara}(i).
\begin{cor} 
Let $({\mathfrak m},V,\omega,B_{\mu})$ be a real   SSR such that $Q\not\equiv0$. Then $({\mathcal F},{\mathcal G})$ is a  one-dimensional abstract noncommutative, completely integrable system on ${\mathcal O}$.
\end{cor}
\begin{demo} Since ${\mathcal O}\subseteq V$ is open and not empty, $({\mathcal O},\omega)$ is a symplectic manifold. According to \cite{LMV} Definition 5.1, $ ({\mathcal F},{\mathcal G})$ is a  one-dimensional abstract noncommutative, completely integrable system on $({\mathcal O},\omega)$  iff
\begin{itemize}
\item[(1)]   ${\mathcal F}$ is of rank $1$ and   ${\mathcal G}$ is of corank $1$.

\item[(2)] ${\mathcal F}\subseteq{\mathcal G}$.

\item[(3)]  $\omega({\mathcal F},{\mathcal G})=0$.
\end{itemize}
\noindent  By definition, the distribution tangent to ${\mathcal F}$ is $A\mapsto{\rm Ker}\,d\mu_A$ and the distribution tangent to ${\mathcal G}$ is $A\mapsto{\mathfrak m}\cdot A$. Hence the above three conditions are equivalent to:
\begin{itemize}
\item[(1')] For all $A\in{\mathcal O}$, \quad dim$({\rm Ker}\,d\mu_A)=1$  and   codim$({\mathfrak m}\cdot A)=1$.

\item[(2')] For all $A\in{\mathcal O}$, \quad${\rm Ker}\,d\mu_A\subseteq{\mathfrak m}\cdot A$.

\item[(3')] For all $A\in{\mathcal O}$, \quad$\omega({\rm Ker}\,d\mu_A,{\mathfrak m}\cdot A)=0$.
\end{itemize}
\noindent  (1') follows immediately from Proposition \ref{skychara}(i)). (2') and (3') follow from the fact that
${\mathfrak m}\cdot A={\mathfrak m}_{\mu}\cdot A=({\rm Ker}\,d\mu_A)^{\perp}$ (Propositions \ref{skychara}(i) and  \ref{smu}).

\end{demo}
\section{Lagrangian decomposition}

In the 19th and early 20th centuries it was discovered that the set of binary cubics  has many remarkable properties and  in this chapter we show that some of them can be generalised to any other SSR. For instance, we determine explicitly the fibres of the quadratic covariant restricted to generic points (cf Corollary \ref{mufibres1}) and show that the symplectic covariants of an SSR satisfy an  Eisenstein syzygy (cf Theorem \ref{Eisen2}). The main results are   Theorem \ref{sum1} and Theorem \ref{quadextsum1} which generalise to all SSRs  decomposition theorems proved  by L.E Dickson for complex binary cubics
and by N.J. Hitchin for real or complex alternating $3$-forms in six dimensions.
\subsection{Symplectic relations of an $SSR$}

Hermite  \cite{He} investigated integral binary cubic forms having the same Hessian. Later G.B. Mathews \cite{Ma} proved some identities for integral binary cubics that essentially determined the projective fibre of the determinant of the Hessian (cf. Example \ref{binarycubics}).  In this section we will obtain generalisations of all these identities in the context of SSRs over $k$. 

\begin{prop}\label{invs1}
Let $({\mathfrak m},V,\omega,B_{\mu})$ be an SSR and let
$A\in V$. Then

\medskip
\noindent(a) $\mu \left(\Psi(A)\right)=-Q(A)\mu(A)$.

\medskip
\noindent(b) $\Psi\left( \Psi(A)\right)=
-Q(A)^2A$.

\medskip
\noindent(c) $Q(\Psi(A))=Q(A)^3$.
\end{prop}
\begin{demo} By definition,
$$
\mu \left(\Psi(A)\right)=\mu \left( \mu(A)\cdot A\right)=B_{\mu}\left(  \mu(A)\cdot A, \mu(A)\cdot A\right)
$$
and by equivariance of the moment map,
$$
B_{\mu}\left( \mu(A)\cdot A, \mu(A)\cdot A\right)=
\left[\mu(A), B_{\mu}\left(A, \mu(A)\cdot A\right)\right]
-B_{\mu}\left( A,\mu(A)\cdot (\mu(A)\cdot A)\right).
$$
The first term vanishes since
\begin{equation}\label{Apsiortho}
B_{\mu}\left(A,\Psi(A)\right)=B_{\mu}\left(A, \mu(A)\cdot A\right)=\frac{1}{2}[\mu(A),\mu(A)]=0
\end{equation}
 and hence
\begin{equation}\label{mompsi}
\mu \left(\Psi(A)\right)=-B_{\mu}\left( A,\mu(A)\cdot (\mu(A)\cdot A)\right).
\end{equation}
To simplify the RHS we first write
$$
\mu(A)\cdot (\mu(A)\cdot A)=B_{\mu}(A,A)\cdot(\mu(A)\cdot A)
$$
and then use  equations \eqref{ssrdefi} and \eqref{Apsiortho}   to get 
\begin{align}\label{mompsi1}
\nonumber
B_{\mu}(A,A)\cdot(\mu(A)\cdot A)&=
B_{\mu}(A,\mu(A)\cdot A)\cdot A
+\frac{3}{ 2}\omega(A,\mu(A)\cdot A)A\\
&=\frac{3}{2}\omega(A,\mu(A)\cdot A)A
=Q(A)A.
\end{align}
Substituting into \eqref{mompsi} above gives (a).

To prove (b) we have
\begin{align}
\nonumber
\Psi\left( \Psi(A)\right)&=
\mu(\mu(A)\cdot A)\cdot\left(\mu(A)\cdot A\right)\\
\nonumber
&=-Q(A)\mu(A)\cdot\left(\mu(A)\cdot A\right)
\quad(\text{by (a)})\\
\nonumber
&=-Q(A)^2A
\quad(\text{by \eqref{mompsi1}}).
\end{align}
Part (c) follows from $Q(\Psi(A))=\frac{3}{2}\omega(\Psi(A),\Psi^2(A))$
and (b).
\end{demo}

As a corollary, we  calculate the values of the three covariants on linear combinations of $A$ and $\Psi(A)$ for $A\in V$. 

\begin{cor}\label{val3fundinv}
Let $({\mathfrak m},V,\omega,B_{\mu})$ be an SSR and let
$A\in V$.

\medskip
\noindent(i) For all $a,a',b,b'\in k$, 
$$
B_{\mu}(aA+b\Psi(A),a'A+b'\Psi(A))=(aa'-Q(A)bb')\mu(A).
$$
\noindent(ii) For all $a,b,\in k$, 
$$\Psi(aA+b\Psi(A))=
(a^2-Q(A)b^2)\left( Q(A)bA+a\Psi(A)\right).
$$
\noindent(iii) For all $a,b,\in k$, \,$Q(aA+b\Psi(A))=(a^2-Q(A)b^2)^2Q(A)$.

\medskip
\noindent(iv) For all $a,b,\in k$, \,$\mu(A)\cdot(aA+b\Psi(A))=a\Psi(A)+Q(A)bA$.

\medskip
\noindent(v) For all $X\in<A,\Psi(A)>$,\quad$\mu(A)\cdot\left(\mu(A)\cdot(X)\right)=Q(A)X$.

\end{cor}
\begin{demo} 
Part (i)  follows from $B_{\mu}\left(A,\Psi(A)\right)=0$ and (a). Part (ii) follows from (i) and  equation \eqref{mompsi1}. Part (iii) follows from (i),  part (iv) from equation \eqref{mompsi1} and part (v) from equation \eqref{mompsi1}.
\end{demo}

\subsection{Lagrangian decomposition}

We come to the key  property of SSRs. Let $P$ be a homogeneous polynomial of degree $3$ in two variables  over an algebraically closed field of characteristic zero.
It was known to  L.E Dickson \cite{Di} p. 17 that  if the discriminant of $P$  is nonzero,  there is a decomposition
$$
P=P_1+P_2
$$
where $P_1$ and $P_2$ are cubes of linear forms. N.J. Hitchin  \cite{Nigel} showed that a suitably generic real or complex exterior $3$-form in six dimensions $\rho$ has a decomposition
$$
\rho=\rho_1+\rho_2
$$
where $\rho_1$ and $\rho_2$ are {\it decomposable} $3$-forms and that this decomposition is unique up to permutation of $\rho_1$ and $\rho_2$. In this section we prove a decomposition theorem which generalises these two results to an arbitrary SSR; moreover, the proof gives  an explicit formula for the summands in terms of  the symplectic covariants.

\begin{lem} \label{Psiofsum} Let  $({\mathfrak m},V,\omega,B_{\mu})$  be an SSR and let $B,C\in V$ be such that $\mu(B)=\mu(C)=0$. Then

\medskip
\noindent(i) $\Psi(B+C)=3\omega(B,C)(-B+C)$;

\medskip
\noindent(ii) $ Q(B+C)=(3\omega(B,C))^2$.
\end{lem}
\begin{demo}  To prove part (i), we have:
\begin{align}
\nonumber
\Psi(B+C)&=\mu(B+C)\cdot (B+C)\\
\nonumber
&=B_{\mu}(B+C,B+C)\cdot (B+C)\\
\nonumber
&=2B_{\mu}(B,C)\cdot B+2B_{\mu}(B,C)\cdot C\quad\text{(since }\mu(B)=\mu(C)=0)\\
\nonumber
&=2\omega(C,B)B-\omega(B,C)B+2\omega(B,C)C-\omega(C,B)C\text{ (by  \eqref{ssrdefi})}\\
\nonumber
&=3\omega(B,C)(-B+C).
\end{align}
As for part (ii), by definition,
$$
Q(B+C)=\frac{3}{2}\omega(B+C,\Psi(B+C))
$$
and by  (i) this gives
$$
Q(B+C)=\frac{9}{2}\omega(B,C)\omega(B+C,-B+C) = 9\omega(B,C)^2.
$$
\end{demo}

\begin{theo}\label{sum1}
Let  $({\mathfrak m},V,\omega,B_{\mu})$  be an SSR. For $A\in V$. The following are equivalent:

\medskip
\noindent(i)   $Q(A)\in{k^*}^2$. 

\medskip
\noindent(ii)  There exist $B,C\in V$ such that
$\mu(B)=\mu(C)=0$, $\omega(B,C)\not=0$ and $A=B+C$.

\medskip
\noindent Moreover, when (i) holds $B$ and $C$ of (ii) satisfy $(3\omega(B,C))^2=Q(A)$,  and  there is a 
square root $q\in k^*$ of $Q(A)$ such that
$$
B=\frac{1}{2}(A+\frac{1}{q}\Psi(A)),\quad C=\frac{1}{2}(A-\frac{1}{q}\Psi(A)).
$$
In particular, $B$ and $C$ are  unique up to permutation.
\end{theo}
\begin{demo}
$(i)\Rightarrow(ii)$: Choose  $q\in k^*$ such that $Q(A)=q^2$.
Set $B=\frac{1}{2}(A+\frac{1}{q}\Psi(A))$ and $C=\frac{1}{2}(A-\frac{1}{q}\Psi(A))$. Then
$\mu(B)=\mu(C)=0$ by Corollary \ref{val3fundinv}(i) and 
$$
\omega(B,C)=-\frac{1}{2q}\omega(A,\Psi(A))=-\frac{1}{3q}Q(A)=-\frac{q}{3}\not=0.
$$

\noindent $(ii)\Rightarrow(i)$:  Suppose now  $A=B+C$ where $B,C\in V$ satisfy $\mu(B)=\mu(C)=0$. It is immediate from  Lemma  \ref{Psiofsum} that
\begin{equation}\label{sqrootQ}
Q(A)=(3\omega(B,C))^2
\end{equation}
which proves (i).

To prove that $B$ and $C$ in the decomposition (ii) are unique up to permutation, suppose $A=B'+C'$ with $B',C'$ satisfying the 
properties of (ii). 
By  Lemma \ref{Psiofsum}, $\Psi(A)=3\omega(B',C')(-B'+C')$  and so $\{B',C'\}$ is a basis of $<A,\Psi(A)>$. But for $X\in<A,\Psi(A)>$, we have $\mu(X)=0$ iff $X$ is proportional to $A+\frac{1}{q}\Psi(A)$ or to $A-\frac{1}{q}\Psi(A)$ (cf Corollary \ref{val3fundinv}(i)). It follows immediately that $\{B',C'\}=\{\frac{1}{2}(A+\frac{1}{q}\Psi(A)),\frac{1}{2}(A-\frac{1}{q}\Psi(A))\}$. This proves the desired uniqueness.
\end{demo}
\begin{rema}
 In general, one cannot distinguish  $B$ and $C$.
However if the field $k$ has a `square root' map, i.e., if there exists a homomorphism
$^{\sqrt{\quad}}:{k^*}^2\rightarrow k^*$ such that  $(\sqrt{x})^2=x$ for all  $x\in {k^*}^2$,
then there is a unique ordered pair  $(B,C)$ such that 
$$
A=B+C\quad\text{and}\quad\omega(B,C)=\sqrt{\omega(B,C)^2}.
$$ 
Examples of fields with a square root map are  the real numbers and a finite field with $n$ elements if  $n=3 $ (mod 4). The complex numbers do not have a square root map.
\end{rema}

The following corollary shows that  the quartic $Q$ takes either no non-zero square values or all non-zero square values. This should be contrasted to the case of coefficients in, say, $\bb Z$.
\begin{cor} 
Either ${\rm Im\,}Q\cap {k^*}^2=\emptyset$ or ${\rm Im\,}Q\cap {k^*}^2= {k^*}^2$.
\end{cor}
\begin{demo} Suppose ${\rm Im\,}Q\cap {k^*}^2\not=\emptyset$.  Then there exists  $A\in V$ and $r\in{k^*}$ such that $Q(A)=r^2$. By the theorem,  there exist $B,C\in V$ such that $\mu(B)=\mu(C)=0$, $A=B+C$  and $r=3\omega(B,C)$. Let $\lambda={r'}^2\in{k^*}^2$  and set
$$
A'=\frac{r'}{r}B+C.
$$
Then by Lemma  \ref{Psiofsum}, $Q(A')=(3\omega(\frac{r'}{r}B,C))^2=\lambda$.
\end{demo}

Theorem \ref{sum1} has an analogue  when $Q(A)$ is not a square in $k^*$ but we need some notation before stating it. We omit the  proof which is straightforward.
Let $k'$ be a quadratic extension of $k$. Since char$(k)\not=2$, the extension $k'/k$ is Galois and  the Galois group ${\rm Gal}(k'/k)$ is isomorphic to $\mathbb Z_2$.  If $W$  is a $k$-vector space, the Galois group acts naturally on $W'=W\otimes_{k}k'$ and we always denote the action of the generator by $w\mapsto \bar{w}$. 
\begin{theo}\label{quadextsum1}
Let  $({\mathfrak m},V,\omega,B_{\mu})$  be an SSR, 
let $\lambda\in k^*\setminus{k^*}^2$, let $k'$ be a splitting field of $x^2-\lambda$ and let
$({\mathfrak m}',V',\omega',B_{\mu}')$  be the SSR obtained by base extension.
For $A\in V$ the following are equivalent:

\medskip
\noindent(i)   $Q(A)\in \lambda{k^*}^2$. 

\medskip
\noindent(ii)  There exist $B\in V'$ such that
$\mu'(B)=0$, $\omega'(B,\overline{B})\not=0$ and 
$A=B+\overline{B}.$

\medskip
\noindent Moreover when (i) holds, $B$  of (ii)  satisfies $(3\omega'(B,\overline{B}))^2=Q(A)$  and there is a  square root $q\in {k'}^*$ of $Q(A)$ such that
$$
B=\frac{1}{2}(A+\frac{1}{q}\Psi(A)).
$$
In particular,  $B$  of (ii) is  unique up to conjugation.
\end{theo}

  \begin{ex}\label{e7sumofpure} For the  SSR given by the half-spinors of a 12-dimensional hyperbolic quadratic form (Example \ref{E7spinconstr})  it follows directly from the definition of $B_{\mu}$ and p108 of \cite{Elie} that the zero set of $\mu$ exactly coincides with the set of pure spinors defined by Cartan.  Thus Theorem \ref{sum1} takes the form: given $A$ a half-spinor  in twelve dimensions with $Q(A)$ a nonzero square in $k^*$,  there are pure spinors $P_1,P_2$ unique up to permutation  such that $A=P_1+P_2$. This seems to be a  fact about spinors which was not known to Cartan.
   \end{ex}
  \begin{ex}\label{e6sumofpure}
  For the 20-dimensional SSR given by  three forms in six dimensions one can check ( see \cite{Nigel}  if $k=\bb R$  or $k=\bb C$ )  that the zero set of $\mu$  is the set of decomposable three forms. Thus Theorem \ref{sum1} takes the form: given  $A$ a three form  in 6 dimensions  with $Q(A)$ a nonzero square in $k^*$,  there are decomposable three forms  $P_1,P_2$ unique up to permutation such that  $A=P_1+P_2$. 
   \end{ex}
     \begin{ex}\label{f4sumofpure}
 For the 14-dimensional SSR given by the primitive three forms of a   symplectic 6-dimensional vector space it follows from the previous example  that the zero set of $\mu$  is the set  of decomposable three forms which are Lagrangian, i.e., whose annihilator is  a Lagrangian subspace. Thus Theorem \ref{sum1} takes the form: given $A$ a primitive three form   in 6 dimensions with $Q(A)$ a nonzero square in $k^*$,  there are decomposable Lagrangian three forms  $P_1,P_2$ unique up to permutation such that  $A=P_1+P_2$.
   \end{ex}
\begin{ex}\label{g2sumofpure}
For the 4-dimensional SSR given by homogeneous polynomials of degree three on a 2-dimensional vector space it is shown in (see \cite{Sl-St1}) 
that the zero set of the moment map  is given by
$$
\mu(P)=0\quad\Leftrightarrow\quad\exists \lambda\in k,\,\alpha\in{k^2}^*\text{ such that }
P=\lambda\alpha^3.
$$
Theorem \ref{sum1} then says that a binary cubic $P$ can be written as a linear combination  of the cubes of two independent linear forms
$$
P=\lambda\alpha^3+\lambda'\alpha'^3
$$
iff $Q(P)$  is a nonzero square, and then $\lambda\alpha^3$ and $\lambda'\alpha'^3$ are unique up to permutation.  This generalises a result of L.E Dickson who showed in \cite{Di} that a complex binary cubic of
nonzero discriminant can be written as the sum of two cubes of linear forms. In \cite{Sl-St1} the authors show that the map
$$
P\mapsto [(\omega(\lambda\alpha^3,\lambda'\alpha'^3),[\lambda\lambda'^{-1}]]\in k^*\times_{\bb Z_2}k^*/{k^*}^3
$$
factors to define a bijection from the set of $SL(2,k)$-orbits of binary cubics for which $Q$ is a nonzero square  to   $k^*\times_{\bb Z_2}k^*/{k^*}^3$.
  \end{ex}

\subsection {Explicit description of the  fibre of $\mu$}

An important consequence of Theorem \ref{sum1} in the finite dimensional case is that if  $Q(P)$ is a nonzero square then $\mu(P)$ is diagonalisable on $V$. We need the following notation (here $\lambda\in k^*$):
 
 $$
\begin{array}{rlc}
{\mathcal O}&:=\{A\in V: &Q(A)\not=0\},\\
{\mathcal O}_{\lambda}&:=\{A\in V: &Q(A)\in\lambda{k^*}^2 \},\\
Z&:=\{A\in V\setminus\{0\}: &\mu(A)=0 \}.
\end{array}
$$ 

\begin{theo} \label{mudiag}
Let  $({\mathfrak m},V,\omega,B_{\mu})$  be an SSR.
Let $A\in{{\mathcal O}}_{1}$ and let $A=B+C$ with $B,C\in Z$ as in  Theorem \ref{sum1}.
\vskip0.1cm
\noindent(i) The restriction of $\mu(A)$ to the subspaces of  $V$
\begin{equation}\label{diagmu}
<B>, \quad C^{\perp}\cap{\rm Ker}\,d\mu_B,\quad B^{\perp}\cap{\rm Ker}\,d\mu_C,\quad<C>
\end{equation}
acts   respectively by the scalars
$$
-3\omega(B,C),\quad-\omega(B,C),\quad\omega(B,C),\quad3\omega(B,C).
$$
\vskip0.1cm
\noindent(ii) There is a direct sum decomposition
$$
V=<B> \oplus \,C^{\perp}\cap{\rm Ker}\,d\mu_B\oplus \,B^{\perp}\cap{\rm Ker}\,d\mu_C \,\oplus <C>.
$$
\vskip0.1cm
\end{theo}
\begin{demo} 
Let $X\in {\rm Ker}\,d\mu_B$.   Equation \eqref{ssrdefi}
$$
2B_{\mu}(B,C)\cdot X-2B_{\mu}(B,X)\cdot C=
2\omega(C,X)B-\omega(B,C)X+\omega(B,X)C
$$
reduces to
$$
\mu(A)\cdot X=
2\omega(C,X)B-\omega(B,C)X
$$
since 
$$
\mu(A)=B_{\mu}(B+C,B+C)=\mu(B)+2B_{\mu}(B,C)+\mu(C)=2B_{\mu}(B,C)
$$
 and $ {\rm Ker}\,d\mu_B$ is isotropic (cf Proposition \ref{skychara}(ii)). Taking  $X=B$ and then $X\in C^{\perp}\cap{\rm Ker}\,d\mu_B$,  this gives
$$
\mu(A)\cdot B=-3\omega(B,C)B,\quad\mu(A)\cdot X=-\omega(B,C)X
$$
respectively. Similarly,  we can show that
$$
\mu(A)\cdot C=3\omega(B,C)C,\quad \mu(A)\cdot Y=\omega(B,C)Y
$$
if $Y\in B^{\perp}\cap{\rm Ker}\,d\mu_C$ and this proves (i).

To prove (ii) we first remark that the   sum
$$
<B>+\, C^{\perp}\cap{\rm Ker}\,d\mu_B\,+\, B^{\perp}\cap{\rm Ker}\,d\mu_C\,+<C>
$$
is direct since the summands correspond to distinct eigenvalues of $\mu(A)$. Hence to prove that this sum is equal to $V$, we  have to prove that the sum of the dimensions  of the summands is equal to $2n$, the dimension of $V$.  Since 
$$
{\rm dim\,}(C^{\perp}\cap{\rm Ker}\,d\mu_B)\ge{\rm dim\,}C^{\perp}+{\rm dim\,Ker}d\mu_B-{\rm dim\,}V
$$
and ${\rm Ker}\,d\mu_B$ is Lagrangian (cf Proposition \ref{skychara}(ii)), we have
$$
n\ge{\rm dim\,}(C^{\perp}\cap{\rm Ker}\,d\mu_B)\ge 2n-1+n-2n=n-1.
$$
However $B\in{\rm Ker}\,d\mu_B\setminus C^{\perp}$ so
${\rm dim\,}(C^{\perp}\cap{\rm Ker}\,d\mu_B)=n-1$.  Similarly one shows that ${\rm dim\,}(B^{\perp}\cap{\rm Ker}\,d\mu_C)=n-1$ and (ii) is proved.
\end{demo}

\begin{rema}\label{sqrteigspaces}
 If $q$ is a square root  of $Q(A)$,   the eigenvalues of  $\mu(A)$  are: $-q,-\frac{q}{3},\frac{q}{3},q$. The sums of eigenspaces   $E_{-q}\oplus E_q$ and $E_{-\frac{q}{3}}\oplus E_{\frac{q}{3}}$ are independent of the choice of square root and in fact by Theorem \ref{mudiag}  and Lemma \ref{Psiofsum},
$$
E_{-q}\oplus E_q= <A,\Psi(A)>, \quad E_{-\frac{q}{3}}\oplus E_{\frac{q}{3}}=<A,\Psi(A)>^{\perp}.
$$ 
\end{rema}
We can now give an explicit description of the fibres of $\mu:{{\mathcal O}}_{\lambda}\rightarrow \mathfrak{m}$.
\begin{cor} \label{mufibres1}
(i) Let $A\in{{\mathcal O}}$. Then
$$
\mu^{-1}(\mu(A))=
\{xA+y\Psi(A):\, x^2-Q(A)y^2=1\}.
$$
\vskip0.1cm
\noindent(ii)  If $A\in{{\mathcal O}}_{1}$ and $A=B+C$ with $B,C\in Z$ as in Theorem \ref{sum1}, then
$$
\mu^{-1}(\mu(A))=\{u B+\frac{1}{u} C:\,u\in k^*\}.
$$
\vskip0.1cm
\noindent(iii) If $A\in{{\mathcal O}}_{\lambda}$ where $\lambda\in k^*\setminus{k^*}^2$ and $A=B+\overline{B}$ with $B\in Z'$ as in Theorem \ref{quadextsum1}, then
$$
\mu^{-1}(\mu(A))=\{z B+\frac{1}{z} \overline{B}:\,z\in k' \text{ such that }z\bar{z}=1\}.
$$
\end{cor}
\begin{demo} First let us remark that parts (ii) and (iii) follow from  part (i) since $\Psi(A)=3\omega(B,C)(-B+C)$ in the first case and $\Psi(A)=3\omega'(B,\overline{B})(-B+\overline{B})$ in the second case.

To prove part (i) we can clearly assume that $Q(A)$ is a square in $k^*$.
 Let $A'\in\mu^{-1}(\mu(A))$. Then $\mu(A')=\mu(A)$ implies
 $Q(A')=Q(A)$ and so $A'\in{{\mathcal O}_1}$.  Let $q$ be a square root of $Q(A')=Q(A)$ and let $E_{\pm q}$ be the eigenspaces of $\mu(A')=\mu(A)$ corresponding to the eigenvalues $\pm q$. By Remark  \ref{sqrteigspaces},
 $$
 E_q\oplus E_{-q}=<A',\Psi(A')>=<A,\Psi(A)>
 $$
 and hence $A'=xA+y\Psi(A)$ for some $x,y, \in k$. Finally, from $\mu(A')=\mu(A)$ we get $x^2-Q(A)y^2=1$ (cf Proposition \ref{val3fundinv}).
 \end{demo}
 \begin{cor} \label{musemisimp}
Let  $({\mathfrak m},V,\omega,B_{\mu})$  be a special symplectic representation and let $A\in V$.  
\vskip.1cm
\noindent(i) If $Q(A)\not=0$ the minimal polynomial of $\mu(A)$ acting on $V$  is $(x^2-Q(A))(x^2-\frac{1}{9}Q(A))$. In particular,
$\mu(A)$ is diagonalisable after at most quadratic extension.
\vskip.2cm
\noindent(ii) If $Q(A)=0$ then  $\mu(A)^4=0$.
\end{cor}
\begin{demo} Part (i) follows immediately from the theorem  and so does the fact that the affine variety
$$
X_{(k)}=\{A\in V:\quad  (\mu(A)^2-Q(A)Id)(\mu(A)^2-\frac{1}{9}Q(A)Id)=0 \}
$$
contains the Zariski open set $Z_{(k)}=\{A\in V:\, Q(A)\not=0\}$. This is also true for the SSR $({\mathfrak m}\otimes_{k}\bar{k},V\otimes_{k}\bar{k},\omega\otimes_{k}\bar{k},B_{\mu}\otimes_{k}\bar{k})$ where
$\bar{k}$  is the algebraic closure of $k$ so $Z_{(\bar{k})}\subseteq X_{(\bar{k})}$. If $Q\not\equiv 0$ this implies  $X_{(\bar{k})}=V\otimes_{k}\bar{k}$ and hence $X_{(k)}=X_{(\bar{k})}\cap V=V$. This means 
$$
 (\mu(A)^2-Q(A)Id)(\mu(A)^2-\frac{1}{9}Q(A)Id)=0\quad\forall A\in V
$$
and in particular $\mu(A)^4=0$ if $Q(A)=0$.  If  $Q\equiv 0$ , $\mu(A)=\tau(A)$ (cf Proposition \ref{quadcovvan}) so $\mu(A)^2=0$.
\end{demo}
 \subsection{Generalized  Eisenstein syzygy for  an $SSR$}
 Consider the equation
 \begin{equation}\label{Eisen1}
x^2-\Delta y^2 = 4z^3.
 \end{equation}
G.  Eisenstein showed in \cite{E} that the values of an integral  binary cubic and its three (classical) covariants at any point $v\in  \bb Z^2$ provide a solution $(x,y,z,\Delta)$ of the equation and Mordell \cite{Mo} p. 216 proved essentially the converse.
 
 In \cite{Sl-St1} we gave a  formulation  of this relation  expressed only in terms of the symplectic covariants of the space of  binary cubics  viewed as an SSR. In this section we prove an identity   satisfied by the symplectic covariants of  any  SSR   and show how it generalises the Eisenstein identity for binary cubics. This is to be contrasted to the statement in \cite{Kab} p. 4657 : \lq  there can
be no analogue of the syzygy (1.1) in general \rq. 

\begin{theo} \label{Eisen2}
Let $({\mathfrak m},V,\omega, B_{\mu})$ be  an SSR. For all $P\in V$, the following identity holds in $ {\mathfrak {sp}}(V,\omega)$ :
\begin{equation}\label{genEis}
\tau(\Psi(P))-Q(P)\tau(P)=-\frac{3}{4}\mu(P)^3+ \frac{1}{12}Q(P)\mu(P).
\end{equation}
\end{theo}
\begin{demo} If  $Q$ vanishes identically,   $\Psi$ vanishes identically,
$\mu(P)=\tau(P)$ (Proposition \ref{quadcovvan}) and  $\mu(P)^2=0$. All terms in \eqref{genEis} vanish and the identity is true.

Suppose $Q$ does not vanish identically. We can suppose without loss of generality that $k$ is algebraically closed and our strategy will be first to prove the identity for $P$ such that $Q(P)\ne 0$, and then to deduce the general case by Zariski closure.

Fix $P\in  V$ such that $Q(P)\ne 0$. Then $P$ and $\Psi(P)$ are linearly independent and we have the symplectic orthogonal decomposition
\begin{equation}\label{perpdecomp}
 V = <P,\Psi(P)> \oplus <P,\Psi(P)>^{\perp}
\end{equation}
which is stable under the action of  $\tau(P)$, $\tau(\Psi(P))$ and $\mu(P)$. Hence to prove  \eqref{genEis} it is sufficient to evaluate on vectors which are either in  $<P,\Psi(P)> $ or in $<P,\Psi(P)>^{\perp}$.

If $X\in<P,\Psi(P)>^{\perp}$, both terms of the LHS  of  \eqref{genEis} evaluated on $X$  give zero. The RHS on $X$ also gives zero since $\mu(P)^2=\frac{Q(P)}{9}Id$ on $<P,\Psi(P)>^{\perp}$ by Remark \ref{sqrteigspaces}.

If $X=aP+b\Psi(P)$,   then  
$$\tau(\Psi(P))(X)=a\omega(\Psi(P),P)\Psi(P), \quad \tau(P)(X)=b\omega(P,\Psi(P))P
$$ 
so the LHS of \eqref{genEis} on $X$ gives
$$
\omega(\Psi(P),P)(a\Psi(P)+bQ(P)P)=-\frac{2}{3}Q(P)\mu(P)
$$
using $\omega(\Psi(P),P)=-\frac{2}{3}Q(P)$ (cf definition) and $\mu(P)=a\Psi(P)+bQ(P)P$ (cf Proposition \ref{val3fundinv}). The RHS of \eqref{genEis}
on $X$ is
$$
(-\frac{3}{4}+\frac{1}{12})Q(P)\mu(P)=-\frac{2}{3}Q(P)\mu(P)
$$
since $\mu(P)^2=Q(P){\rm Id}$ on $<P,\Psi(P)>$ by Remark \ref{sqrteigspaces}.

We have now proved (\ref{genEis}) in the case where $P$ satisfies $Q(P)\ne 0$. The difference between the LHS and the RHS of (\ref{genEis}) therefore defines a polynomial function of degree six on $V$, say $\Delta$, which vanishes on the non-empty Zariski open set ${\mathcal O} = \{P\in V:Q(P)\ne 0\}$. By continuity, $\Delta$ vanishes on the Zariski closure of ${\mathcal O}$ which is $V$ since non-empty Zariski open sets are dense if $k$ is algebraically closed. This completes the proof of (\ref{genEis}).
\end{demo}

We saw in  Corollary \ref{musemisimp} that  $\mu(P)$ acting on $V$ is nonzero semisimple if $Q(P)\not=0$ and that $\mu(P)^4=0$ if  $Q(P)=0$. In the second case we can now say something about the nilpotency index of $\mu(P)$
\begin{cor} 
Let $({\mathfrak m},V,\omega, B_{\mu})$ be  a special symplectic representation and let $P\in V$.
Then $\mu(P)^3=0$ iff $\Psi(P)=0$.
\end{cor}
\begin{demo} If $\mu(P)^3=0$ we must have $Q(P)=0$ since otherwise, $\mu(P)$ would be nonzero semisimple by Corollary \ref{musemisimp}. Equation \eqref{genEis} then reduces to
$$
\tau(\Psi(P))=0
$$
which implies $\Psi(P)=0$ since $\tau:V\rightarrow sp(V,\omega)$ is injective.

To prove implication in the other direction, suppose $\Psi(P)=0$. Then $Q(P)=\frac{3}{2}\omega(P,\Psi(P))=0$ and   \eqref{genEis}  reduces to $\mu(P)^3=0$.
\end{demo}

 In the case of the special symplectic representation corresponding to binary cubics (cf Example \ref{binarycubics}), we now show that the identity \eqref{genEis} implies the classical Eisenstein identity.
 
  Let $(V,\Omega)$ be a two-dimensional symplectic vector space, let $\,{\tilde{}}:V\mapsto V^*$ be the isomorphism  defined by ${\tilde{v}}(w)=\Omega(v,w)$ and give  $S^3(V^*)$, the space of cubic functions  on $V$, the  unique symplectic structure $\omega$ satisfying
\begin{equation}\label{symp+ev}
P(v)=\omega (P, {\tilde{v}}^3)\quad\forall v\in V,\,\forall P\in S^3(V^*).
\end{equation}
Evaluating  \eqref{genEis}  at ${\tilde{v}}^3$ and contracting with ${\tilde{v}}^3$:
\begin{equation}\label{Eisred1}
\begin{array}{rl}
\omega(\tau(\Psi(P))\cdot {\tilde{v}}^3,{\tilde{v}}^3)-&
Q(P)\omega(\tau(P)\cdot{\tilde{v}}^3,{\tilde{v}}^3)=\\
&-\frac{3}{4}\omega(\mu(P)^3\cdot{\tilde{v}}^3,{\tilde{v}}^3)
+ \frac{1}{12}Q(P)\omega(\mu(P)\cdot{\tilde{v}}^3,{\tilde{v}}^3).
\end{array}
\end{equation}
The LHS is 
$$
\omega(\omega(\Psi(P),{\tilde{v}}^3)\Psi(P),{\tilde{v}}^3)-Q(P)\omega(\omega(P,{\tilde{v}}^3)P,{\tilde{v}}^3)
$$
which by \eqref{symp+ev} simplifies to
\begin{equation}\label{EisexLHS}
\Psi(P)(v)^2-Q(P)P(v)^2.
\end{equation}
To calculate the RHS, we  use the derivation rule   to get
$$
\begin{array}{rl}
\mu(P)\cdot{\tilde{v}}^3&=3(\mu(P)\cdot{\tilde{v}}) {\tilde{v}}^2,\\
\mu(P)^3\cdot{\tilde{v}}^3&=3(\mu(P)^3\cdot\tilde{v}){\tilde{v}}^2+
18( \mu(P)^2\cdot\tilde{v})(\mu(P)\cdot\tilde{v}) {\tilde{v}}+
6(\mu(P)\cdot\tilde{v}) ^3
\end{array}
$$
from which it  follows that
$$
\omega(\mu(P)\cdot{\tilde{v}}^3,{\tilde{v}}^3)=0,\quad
\omega(\mu(P)^3\cdot{\tilde{v}}^3,{\tilde{v}}^3)=6\omega((\mu(P)\cdot\tilde{v}) ^3,{\tilde{v}}^3)
$$
since, by \eqref{symp+ev}, $\omega(R,{\tilde{v}}^3)=0$ if ${\tilde{v}}$ divides $R$. Hence the RHS of \eqref{Eisred1} reduces to
$$
-\frac{9}{2}\omega((\mu(P)\cdot\tilde{v}) ^3,{\tilde{v}}^3)
$$
which by \eqref{symp+ev}, is the cube of the value of the linear form $\mu(P)\cdot\tilde{v}$ at $v$
 (up to a constant):
\begin{equation}\label{EisexRHS}
-\frac{9}{2}\mu(P)(\tilde{v}) ^3(v)=-\frac{9}{2}(\,\mu(P)\cdot\tilde{v}(v)\,)^3.
\end{equation}
Finally, equating  \eqref{EisexLHS} and \eqref{EisexRHS}, we have
$$
\Psi(P)(v)^2-Q(P)P(v)^2=-\frac{9}{2}\,(\mu(P)\cdot\tilde{v}(v))^3\quad\forall v\in V,\forall P\in S^3(V^*).
$$
This is the classical Eisenstein syzygy \eqref{Eisen1}  satisfied by  the values of $P$ and its three covariants at any $v\in V$ since it can be written
$$
x^2-\Delta y^2=4z^3
$$
if we set 
$
y=P(v),\, x=\frac{1}{3}\Psi(P)(v),\, z=-\frac{1}{2}\mu(P)(\tilde{v})(v)$ and $\Delta=\frac{1}{9}Q(P).
$

 \section{Global Lagrangian decomposition and special geometry }
Let  $({\mathfrak m},V,\omega,B_{\mu})$ be an SSR and  let $\lambda\in k^*$.
In this section we  show that     ${{\mathcal O}}_{\lambda}$ carries a `local  geometric structure' which  if $k=\bb R$  reduces to  either a conic, special bi-Lagrangian structure (cf \cite{CLS})   or a conic, special pseudo-K\"ahler structure  (cf \cite{Nigel}, \cite{F}, \cite{ACD}) depending on whether $\lambda$ is a square  or not.  

The main ingredient of  this `local  geometric structure', corresponding to an integrable complex structure  and holomorphic $\bb C^*$-action if $k=\bb R$ and $\lambda=-1$, is obtained by Lagrangian decomposition.
The basic idea is to associate to $P\in{{\mathcal O}}_{\lambda}$  one of the summands given by  Theorem \ref{sum1} or Theorem \ref{quadextsum1}.  Since the summands are in general indistinguishable  we have to go a double cover $p_{\lambda}:\widehat{  {{\mathcal O}}_{\lambda}} \rightarrow{{\mathcal O}}_{\lambda}$ to be able to do this. We then have a map $\alpha:\widehat{  {{\mathcal O}}_{\lambda}}\rightarrow V\otimes_k A_{\lambda}$ which by Lagrangian decomposition takes its values in
$Z_{\lambda}^{gen}$, a Zariski $k$-open  set of solutions  of a system of homogeneous quadratic equations  defined over a quadratic extension $A_{\lambda}$ of $k$. Both $\widehat{  {{\mathcal O}}_{\lambda}}$ and $Z_{\lambda}^{gen}$ are naturally  conic quasi-affine $k$-varieties and we show that $\alpha$ is an isomorphism in this category.  We think of $(\widehat{  {{\mathcal O}}_{\lambda}},\alpha)$ as a chart for a  local geometric structure on $ {{\mathcal O}}_{\lambda}$.

\subsection{Quadratic extensions of $k$ and base extension}
It will be  convenient to have   a uniform description of all quadratic extensions of $k$ which  includes the `degenerate quadratic extension' $k\times k$ (a.k.a the double numbers, the paracomplex numbers, the split-complex numbers, algebraic motors, ...)
\begin{defi}
Let $\lambda\in k^*$. The two-dimensional composition algebra $A_{\lambda}$ is the quotient of the polynomial algebra $k[x]$ by the ideal generated by $x^2-\lambda$:
$$
A_{\lambda}=k[x]/<x^2-\lambda>, 
$$
The conjugation map  \,\,$\bar{}:A_{\lambda}\rightarrow A_{\lambda}$ is induced by $x\mapsto -x$ and the norm of $z\in A_{\lambda}$ is $N(z)=z\bar{z}$. The inclusion $k\hookrightarrow k[x]$ induces a canonical identification of  $k$ with the fixed point set of conjugation  and we set
$ {Im\,A_{\lambda}}=\{z\in A_{\lambda}:\,\bar{z}=-z\}$,
 $A_{\lambda}^*=\{z\in A_{\lambda}: z\text{ is invertible}\}$
 and ${Im^*\,A_{\lambda}}= {Im\,A_{\lambda}}\setminus\{0\}$. We write $\sqrt{\lambda}$ for the class of $x$ in ${Im\,A_{\lambda}}$.
\end{defi}

If $\lambda$ is a square, $A_{\lambda}$ is isomorphic to the  direct product $k\times k$ with conjugation $\overline{(x,y)}=(y,x)$ and norm $N(x,y)=xy$.
If $\lambda$ is not a square, $A_{\lambda}$  is a  splitting field   of $x^2-\lambda$ with
  conjugation given by the action of the Galois group $\bb Z_2$ and $N$ by the norm of the extension $A_{\lambda}/k$. The norm is hyperbolic in the first case and anisotropic in the second.
Finally, note that
  squaring  in $A_{\lambda}$ defines a two to one, surjective map $sq:{Im^*\,A_{\lambda}}\rightarrow \lambda{k^*}^2 $  and that $A_{\lambda}$ is isomorphic to $A_{\lambda'}$ iff $\lambda=\lambda'$ mod ${k^*}^2$.
  
   Let $({\mathfrak m},V,\omega,B_{\mu})$ be an SSR. We
  set $\hat{V}=V\oplus k$ and define $\sigma:\hat{V}\rightarrow\hat{V}$,
${\hat Q}:\hat{V} \rightarrow k$ and  $H:\hat{V} \rightarrow k$    by
$$
\sigma(P,z)=(P,-z),\quad{\hat Q}(P,z)=Q(P), \quad H(P,z)=z.
$$
\begin{defi}
 Let  $\lambda\in k^*$.  Set
$$
\widehat{  {{\mathcal O}}_{\lambda}}=\{\hat{P}\in \hat{V}:{\hat Q}(\hat{P})\not=0 \}\cap
\{\hat{P}\in \hat{V}:{\hat Q}(\hat{P})=\lambda H(\hat{P})^2 \}
$$
\end{defi}
By definition,
$\widehat{  {{\mathcal O}}_{\lambda}}$ is  a quasi-affine $k$-variety, i.e., the intersection of a Zariski open set with a Zariski closed set. Let $S(\hat{V}^*)$ be the ring of polynomial functions  on $\hat{V}$ and let $S(\hat{V}^*)\,[\frac{1}{ {{\hat Q}}}]$  be its localisation at ${{\hat Q}}$ . Restriction to $\widehat{  {{\mathcal O}}_{\lambda}}$ defines a  homomorphism from $S(\hat{V}^*)$ to the ring of functions on $\widehat{  {{\mathcal O}}_{\lambda}}$  and this homomorphism  uniquely extends to a homomorphism defined on $S(\hat{V}^*)\,[\frac{1}{ {{\hat Q}}}]$ since ${{\hat Q}}$ never vanishes on $\widehat{  {{\mathcal O}}_{\lambda}}$. Clearly  
$\frac{ \lambda H}{\hat{Q}}\times H=1$ on  $\widehat{  {{\mathcal O}}_{\lambda}}$ so $\frac{1}{H}$ is a regular function on  $\widehat{  {\mathcal O}_{\lambda} }$ in the sense of the following definition.
\begin{defi} 
The ring $R(\widehat{{\mathcal O}_{\lambda}} )$  of regular functions  on 
$\widehat{  {{\mathcal O}}_{\lambda}}$ is the ring of functions on $\widehat{  {{\mathcal O}}_{\lambda}}$ which are restrictions to $\widehat{  {\mathcal O}_{\lambda} }$ of elements of $S(\hat{V}^*)\,[\frac{1}{ {\hat Q}}]$ .
\end{defi}
There is  obviously a ring isomorphism
$$
S(\hat{V}^*)\,[\frac{1}{\hat { Q}}]\,/\,{\mathfrak I}_{\widehat{  {{\mathcal O}}_{\lambda}}}\cong R(\widehat{{\mathcal O}_{\lambda}} ).
$$
where ${\mathfrak I}_{\widehat{  {{\mathcal O}}_{\lambda}}}$  denotes the ideal of elements of  $S(\hat{V}^*)\,[\frac{1}{ \hat{Q}}]$ which vanish when restricted to $\widehat{  {{\mathcal O}}_{\lambda}}$.
The $k^*$-action on $\hat{V}$ given by
 \begin{equation}\label{wgtedaction}
a\cdot(P,z)=(aP,a^2z)\quad\forall a\in k^*,\forall(P,z)\in \hat{V},
\end{equation}
induces a $\bb Z'$-grading  on  $S(\hat{V}^*)\,[\frac{1}{{ Q}}]$, preserves $\widehat{  {{\mathcal O}}_{\lambda}}$ and  
preserves the ideal  ${\mathfrak I}_{\widehat{  {{\mathcal O}}_{\lambda}}}$ ($\bb Z'=\bb Z$ if $k$ is infinite and  $\bb Z'=\bb Z/(p^n-1)\bb Z$  if  $k$ is finite with $p^n$ elements).
Hence  it also
induces a  $\bb Z'$-grading 
$$ R(\widehat{{\mathcal O}_{\lambda}})=\oplus_{n\in\bb Z'} R_n(\widehat{{\mathcal O}_{\lambda}})
$$ 
where
 $$
 R_n(\widehat{{\mathcal O}_{\lambda}})=\{f\in R(\widehat{{\mathcal O}_{\lambda}}):\, f(a\cdot\hat{P})=a^nf(\hat{P})\,\forall a\in k^*,\forall \hat{P}\in\widehat{{\mathcal O}_{\lambda}}\}.
 $$
We have the following commutative diagram
\begin{equation}\label{hatcover}
\xymatrix
{
\widehat{{\mathcal O}_{\lambda}}\ar[r]^{H{\sqrt{\lambda}}}\ar[d]_{p_{\lambda}}&{\,Im^*\,A_{\lambda}}\ar[d]^{{\rm square}}\\
{\mathcal O}_{\lambda}\ar[r]^{Q}&\lambda{k^*}^2\\
}
\end{equation}
where the projection $p_{\lambda}(P,z)=P$ is two to one surjective and satisfies $p_{\lambda}\circ\sigma=\sigma$. By pullback, $p_{\lambda}$ maps  $S(V^*)[\frac{1}{Q}]$ isomorphically onto the fixed point set  of $\sigma$ acting on $R(\widehat{{\mathcal O}_{\lambda}} )$ and this isomorphism   is compatible with gradings.

We now associate a second conic, quasi-affine $k$-variety to the SSR $({\mathfrak m},V,\omega,B_{\mu})$ and a nonzero scalar $\lambda\in k^*$. As we have seen, $k$ is canonically included in $A_{\lambda}$ so we can base extend any $k$-vector space  $V$ to an $A_{\lambda}$-module $V_{\lambda}=V\otimes_{k}A_{\lambda}$ and any $k$-linear map to an $A_{\lambda}$-morphism.  The action of $\bb Z_2$ by conjugation on $A_{\lambda}$ extends  naturally to an action on  $V_{\lambda}$  whose fixed point set is $V$.
We denote by $({\mathfrak m}_{\lambda},V_{\lambda},\omega_{\lambda},B_{\lambda})$ the base extension of  the special symplectic $k$-representation $({\mathfrak m},V,\omega,B_{\mu})$ to $A_{\lambda}$ and  by ${h}:V_{\lambda}\rightarrow k$ the hermitian quadratic form ${h}(v)=\frac{1}{\sqrt{\lambda}}\omega_{\lambda}(\bar{v},v)$.
\begin{defi} 
Let $\mu_{\lambda}:V_{\lambda}\rightarrow{\mathfrak m}_{\lambda}$ be the quadratic covariant of $({\mathfrak m}_{\lambda},V_{\lambda},\omega_{\lambda},B_{\lambda})$.
Define
$$
Z_{\lambda}^{gen}=\{v\in V_{\lambda}: \,{h}(v)\not= 0\}\cap\{v\in V_{\lambda}:\,  \mu_{\lambda}(v)=0\}.
$$
\end{defi}
Since $\overline{\omega_{\lambda}(\bar{v},v)}=\omega_{\lambda}(v,\bar{v})$ and $\overline{\mu_{\lambda}(v)}=\mu_{\lambda}(\bar{v})$,  conjugation maps $Z_{\lambda}^{gen}$  to $Z_{\lambda}^{gen}$ and has fixed points iff $\lambda$ is a square in $k$.

By definition,
$Z_{\lambda}^{gen}$ is  a quasi-affine $k$-variety, i.e., the intersection of a Zariski open set with a Zariski closed set.
Let $S(V_{\lambda}^*)$ be the ring of $k$-valued $k$-polynomial functions  on $V_{\lambda}$  and let $S(V_{\lambda}^*) [\frac{1}{ {h}}]$ be its localisation at $h$. As in the case of $\widehat{{\mathcal O}_{\lambda}}$ above,  restriction to
$Z_{\lambda}^{gen}$ defines a ring homomorphism from $S(V_{\lambda}^*) [\frac{1}{h}]$ to the ring of functions on $Z_{\lambda}^{gen}$.
\begin{defi} 
The ring of regular functions $R(Z_{\lambda}^{gen} )$  is the ring of functions on 
$Z_{\lambda}^{gen}$ which are restrictions to $Z_{\lambda}^{gen}$ of elements of $S(V_{\lambda}^*) [\frac{1} {h}]$.
\end{defi}

There is evidently a ring isomorphism
$$
S(V_{\lambda}^*) [\frac{1} {h}]\,/\,{\mathfrak I}_{Z_{\lambda}^{gen}}\cong R(Z_{\lambda}^{gen} ).
$$ 
where  ${\mathfrak I}_{Z_{\lambda}^{gen}}$ denotes the ideal of elements  in $S(V_{\lambda}^*) [\frac{1} {h}]$  which vanish when restricted to $Z_{\lambda}^{gen}$.  The natural action of $k^*$ on $V_{\lambda}$ induces a $\bb Z'$-grading on $S(V_{\lambda}^*)[\frac{1}{h}]$ and, since
${\mathfrak I}_{Z_{\lambda}^{gen}}$ is  stable under this action, it also induces
a $\bb Z'$-grading  
 $$ R(Z_{\lambda}^{gen})=\oplus_{n\in\bb Z'} R_N(Z_{\lambda}^{gen})
 $$ 
 where
 $$
 R_N(Z_{\lambda}^{gen})= \{f\in R(Z_{\lambda}^{gen}):\, f(aP)=a^Nf(P)\,\forall a\in k^*,\forall P\in Z_{\lambda}^{gen}\}.
 $$
Similarly,  the ideal ${\mathfrak I}_{Z_{\lambda}^{gen}}\otimes_k  A_{\lambda}\subset S(V_{\lambda}^*)[\frac{1}{h}]\otimes_k A_{\lambda}$  is stable under the natural action of $A_{\lambda}^*$  and this induces  a $(\bb Z'\times\bb Z'$)-grading  
$$R(Z_{\lambda}^{gen})\otimes_k A_{\lambda}=\oplus_{(m,n)\in\bb Z'\times\bb Z'} R_{m,n}(Z_{\lambda}^{gen})
$$ 
where
$$
 R_{m,n}(Z_{\lambda}^{gen})=\{f\in R(Z_{\lambda}^{gen})\otimes_k A_{\lambda}:\, f(aP)=a^m{\bar a}^nf(P)\,\forall a\in A_{\lambda}^*,\forall P\in Z_{\lambda}^{gen}\}.
 $$
 Since the $k^*$-action on $S(V_{\lambda}^*)$ extends the $A_{\lambda}^*$-action on $S(V_{\lambda}^*)\otimes_k A_{\lambda}$, we have a   `type' decomposition:
 $$
  R_N(Z_{\lambda}^{gen})\otimes A_{\lambda}=\oplus_{m+n=N} R_{m,n}(Z_{\lambda}^{gen}).
 $$
 We think of $Z_{\lambda}^{gen}$ as  a   conic (i.e., with a quasi-affine principal $k^*$-action), quasi-affine $k$-variety  but which, because of the way it is defined,  also has  a natural compatible `integrable  $A_{\lambda}$-structure' and compatible `principal holomorphic $A_{\lambda}^*$-action'. For example  if $k=\bb R$ and $ \lambda=-1$, then $Z_{\lambda}^{gen}$ is  a conic real manifold with a natural compatible integrable complex structure (as a real Zariski open set in a complex algebraic variety)  and  principal  holomorphic $\bb C^*$ action.

\subsection{Global Lagrangian decomposition}

To an SSR  $({\mathfrak m},V,\omega,B)$ and  a scalar $\lambda\in k^*$, we have now associated two conic, quasi-affine $k$-varieties: $(\widehat{  {{\mathcal O}}_{\lambda}},R(\widehat{  {{\mathcal O}}_{\lambda}}))$ and $(Z_{\lambda}^{gen},R(Z_{\lambda}^{gen}))$.   
The Lagrangian decomposition theorem provides a natural map from
$\widehat{  {{\mathcal O}}_{\lambda}}$ to $Z_{\lambda}^{gen}$ which we  show is an isomorphism in this category.
\begin{theo}\label{transport}
 Define $\alpha:\widehat{  {{\mathcal O}}_{\lambda}}\rightarrow V_{\lambda}$ 
and $\beta:Z_{\lambda}^{gen}\rightarrow V\oplus k$ by
$$
\begin{array}{rll}
\alpha(P,z)&=\frac{1}{ 2} (P+\frac{1}{z\sqrt{\lambda}}\Psi(P))\qquad&\forall(P,z)\in
\widehat{ {{\mathcal O}}_{\lambda}},\\
\beta (v)&=(v+\bar{v}, \,3{h}(v))\qquad&\forall v \in Z_{\lambda}^{gen}.
\end{array}
$$
(i) $\alpha(\sigma(\hat{P})=\overline{\alpha(\hat{P})}$ for all $\hat{P}\in \widehat{  {{\mathcal O}}_{\lambda}}$.

\noindent(ii) $\alpha$ and $\beta$ commute with the natural $k^*$-actions.

\noindent(iii) $\alpha$ takes  values in $Z_{\lambda}^{gen}$ and $\beta\circ\alpha={\rm Id}_{\widehat{  {{\mathcal O}}_{\lambda}}}$. 

\noindent(iv) $\beta$ takes  values in $\widehat{  {\mathcal O}_{\lambda} }$ and  $\alpha\circ\beta={\rm Id}_{Z_{\lambda}^{gen}}$.

\noindent(v)  $\alpha^*$ maps $R(Z_{\lambda}^{gen})$ to $R(\widehat{  {\mathcal O}_{\lambda} })$.

\noindent(vi) $\beta^*$ maps $R(\widehat{  {\mathcal O}_{\lambda} })$ to $R(Z_{\lambda}^{gen})$.
\end{theo}
\begin{demo} (i) is immediate.

(ii) : 
Let $a\in k$ and $(P,z)\in\widehat{ {{\mathcal O}}_{\lambda}} $. Using equation \eqref{wgtedaction} and the fact that
$\Psi$ is cubic, we have 
$$
\alpha(a\cdot(P,z))=\frac{1}{2}(aP+\frac{1}{a^2z\sqrt{\lambda}}\Psi(aP))=
a\alpha(P,z)
$$
which shows that $\alpha$ commutes with the $k^*$-actions.
It is clear  that $\beta$ commutes with the $k^*$-actions since if $a\in k$ and $v\in V$, $\bar{a}=a$ and $h(av)=a^2v$. This proves (i).

\noindent(iii) : 
 It follows from  Theorems \ref{sum1} and \ref{quadextsum1}  that $\mu_{\lambda}(\alpha(P,z))=0$ since $z\sqrt{\lambda}$ is a square root of $Q(P)$. To prove $\alpha(P,z)\in Z_{\lambda}^{gen}$, we now have  to show that ${h}(\alpha(P,z))\not=0$.
Since ${h}(\alpha(P,z))=\frac{1}{\sqrt{\lambda}}\omega_{\lambda}(\overline{\alpha(P,z)},\alpha(P,z))$ and
\begin{align}
\nonumber
\omega_{\lambda}(\overline{\alpha(P,z)},\alpha(P,z))&=
\frac{1}{4}\left(\omega_{\lambda}(P-\frac{1}{z\sqrt{\lambda}}\Psi(P),P+\frac{1}{z\sqrt{\lambda}}\Psi(P))\right)\\
\nonumber
&=\frac{1}{2 z\sqrt{\lambda}}\omega(P,\Psi(P)),
\end{align}
we have
\begin{equation}\label{pullbackofH}
{h}(\alpha(P,z))=\frac{1}{\sqrt{\lambda}}\times\frac{1}{ 2z\sqrt{\lambda}}\times\frac{2}{3}Q(P)=\frac{z}{3}
\end{equation}
 which is nonzero since $Q(P)=\lambda z^2\not=0$. Hence $\alpha(P,z)\in Z_{\lambda}^{gen}$ and it is now straightforward to check that $\beta(\alpha(P,z))=(P,z)$ so  (ii) is proved.

\noindent(iiv) : 
We first have to show that $\beta$ takes  values in $\widehat{  {\mathcal O}_{\lambda} }$, i.e., that
\begin{equation}\label{betavalues}
{Q}(v+\bar{v})=\lambda\,(3{h}(v))^2.
\end{equation}
The
special symplectic representation $({\mathfrak m}_{\lambda},V_{\lambda},\omega_{\lambda},B_{\lambda})$  was obtained by base extension so if ${Q_{\lambda\,}}$ denotes its quartic covariant, we have
$$
Q(v+\bar{v})={Q_{\lambda\,}}(v+\bar{v})
$$
and since $\mu_{\lambda}(v)=\mu_{\lambda}(\bar{v})=0$,   Lemma \ref{Psiofsum} gives
$$
{Q_{\lambda\,}}(v+\bar{v})=(3\omega_{\lambda}(v,\bar{v}))^2
$$
which proves \eqref{betavalues}. Finally,  to see that $\alpha\circ\beta(v)=v$, we have
$$
\alpha\circ\beta(v)=(v+\bar{v})+\frac{1}{3\omega_{\lambda}(\bar{v},v)}\Psi(v+\bar{v})
$$
and since $\Psi(v+\bar{v})=\Psi_{\lambda}(v+\bar{v})$ and  $\mu_{\lambda}(v)=\mu_{\lambda}(\bar{v})=0$,   by Lemma \ref{Psiofsum} this implies
$$
\alpha\circ\beta(v)=\frac{1}{2}(v+\bar{v})+\frac{1}{3\omega_{\lambda}(\bar{v},v)}(3\omega_{\lambda}(v,\bar{v}))(-\frac{1}{2}v+\frac{1}{2}\bar{v})=v.
$$

\noindent(v) and (vi) :  To prove (iv) it is sufficient to show that
$\alpha^*{h}\in R(\widehat{  {\mathcal O}_{\lambda} })$ and $\alpha^*\eta\in R(\widehat{  {\mathcal O}_{\lambda} })$ for any $\eta\in{\rm Hom}_k(V_{\lambda},k)$ since  ${h}$, $\frac{1}{ h}$ and  restrictions of $k$-linear forms on $V_{\lambda}$ to $Z_{\lambda}^{gen}$ generate $R(Z_{\lambda}^{gen})$.
Similarly, to prove (v) it is sufficient to show that
$\beta^*H\in R(Z_{\lambda}^{gen})$ and $\beta^*\xi\in R(Z_{\lambda}^{gen})$ for any 
$\xi\in V^*$ since  $H$, $\frac{1}{H}$ and  $V^*$ generate $R(\widehat{  {\mathcal O}_{\lambda} })$.
We need the following lemma.
\begin{lem} 
(a) For all $(P,z)\in\widehat{  {{\mathcal O}}_{\lambda}}$ and all $\eta\in{\rm Hom}_k(V_{\lambda},k)$,
\begin{equation}
\alpha^*{h}(P,z)=\frac{1}{3}H(P,z),\qquad\alpha^*\eta(P,z)=\frac{1}{2}\left(\eta(P)+\frac{1}{H(P,z)}\eta(\frac{1}{\sqrt{ \lambda}}\Psi(P))\right).
\end{equation}
(b) For all $v\in Z_{\lambda}^{gen}$ and for all $\xi\in V^*$,
\begin{equation}
\beta^*H(v)=3{h}(v),\qquad \beta^*\xi(v)=\xi( v+\bar{v}  ).
\end{equation}
\end{lem}
\begin{demo} It follows from \eqref{pullbackofH} that
$\alpha^*{h}(P,z)=\frac{z}{3}=\frac{1}{3}H(P,z)$. From the formula for $\alpha$,
$\alpha^*\eta(P,z)=\eta(\alpha(P,z))=\frac{1}{2}\eta(P+\frac{1}{ z\sqrt{\lambda}}\Psi(P))$ and this, since $\eta$ is $k$-linear, simplifies to:
$$
\frac{1}{2}\left(\eta(P)+\frac{1}{ z}\eta(\frac{1}{\sqrt{\lambda}}\Psi(P))\right)
=\frac{1}{2}\left(\eta(P)+\frac{1}{ H(P,z)}\eta(\frac{1}{\sqrt{\lambda}}\Psi(P))\right).
$$
This proves (a) and  (b) is immediate.
\end{demo}

  By the lemma,  $\alpha^*{h}=\frac{1}{3}H$  which is in $R(\widehat{  {\mathcal O}_{\lambda} })$. With the notation  of the lemma,
\begin{equation}\label{regtoreg}
\alpha^*\eta(P,z)=
\frac{1}{2}\left(\eta(P)+\frac{1}{ H(P,z)}\eta(\frac{1}{\sqrt{\lambda}}\Psi(P))\right).
\end{equation}
Since $\eta\in{\rm Hom}_k(V_{\lambda},k)$ and both $V$ and ${1\over\sqrt{\lambda}}V$ are subsets of $V_{\lambda}$ , the functions $(P,z)\mapsto \eta(P)$ and $(P,z)\mapsto \eta({1\over\sqrt{\lambda}}\Psi(P))$ are the restrictions to  $\widehat{  {\mathcal O}_{\lambda} }$ of respectively   linear and  cubic  functions defined on $V\oplus k$. By definition this means they are regular functions on $\widehat{  {\mathcal O}_{\lambda} }$ and since we already know that ${1\over H}$ is regular,  we conclude from \eqref{regtoreg}� that $\alpha^*\eta\in R(\widehat{  {\mathcal O}_{\lambda} })$ and (iv) is proved.

Part(v) is proved similarly. By the lemma $\beta^*H=3{h}$ which is in  $R(Z_{\lambda}^{gen})$. With the notation of the lemma,
$$
 \beta^*\xi(v)=\xi( v+\bar{v}) 
$$
which shows that $\beta^*\xi$ is the restriction to $Z_{\lambda}^{gen}$ of 
$\xi'\in {Hom}_k(V_{\lambda},k)$ defined by $\xi'(x)=\xi(x+\bar{x})$ for $x\in V_{\lambda}$. Hence
$\beta^*\xi\in R(Z_{\lambda}^{gen})$ and (v) is proved.

\end{demo}

Using the maps $\alpha$ and $\beta$, we can transport structure  from $\widehat{  {\mathcal O}_{\lambda} }$ to $Z_{\lambda}^{gen}$ and vice versa.  In particular $Z_{\lambda}^{gen}$  has  a natural  `integrable  $A_{\lambda}$-structure' and  principal  `holomorphic'  $A_{\lambda}^*$-action so  $\widehat{  {\mathcal O}_{\lambda} }$ inherits these structures by transport of structure. 
The principal  $A_{\lambda}^*$-action  on $\widehat{  {\mathcal O}_{\lambda} }$ is given by
$$
a_{\lambda}\cdot\hat{P}=\beta(a_{\lambda}\alpha(\hat{P}))\quad\forall a_{\lambda}\in A_{\lambda}^*,\forall \hat{P}\in\widehat{  {\mathcal O}_{\lambda} },
$$
which explicitly is:
\begin{prop} 
For all $a+b\sqrt{\lambda}\in A_{\lambda}^*$ and for all $(P,z)\in\widehat{  {\mathcal O}_{\lambda} }$,
$$
(a+b\sqrt{\lambda})\cdot(P,z)=\left(aP+{b\over z}\Psi(P),(a^2-b^2\lambda)z\right).
$$
\end{prop}
\begin{demo} This is a  straightforward calculation.
\end{demo}

Using these formulae and Corollary  \ref{mufibres1}, the orbits of the `one-dimensional torus'  group
$$
U(A_{\lambda}^*)=\{a_{\lambda}\in A_{\lambda}^*:\quad a_{\lambda}\overline{a_{\lambda}}=1\}
$$
can be characterised as the level sets of the map $\hat{\mu}:\widehat{  {\mathcal O}_{\lambda} }\rightarrow {\mathfrak m}\times k$ defined by
$$
\hat{\mu}(P,z)=(\mu(P),z)\quad\forall (P,z)\in\widehat{  {\mathcal O}_{\lambda} }.
$$
\begin{cor} 

Let $\hat{P},\hat{P_1}\in \widehat{  {\mathcal O}_{\lambda} }$. There exists
$u\in U(A_{\lambda}^*)$ such that $u\cdot\hat{P}=\hat{P_1}$ iff $\hat{\mu}(\hat{P})=\hat{\mu}(\hat{P_1})$.
\end{cor}
\begin{demo} Let $\hat{P}=(P,z)$ and $\hat{P_1}=(P_1,z_1)$
 If there exists $u\in U(A_{\lambda}^*)$ such that $u\cdot(P,z)=(P_1,z_1)$ then,
writing $u=a+b\sqrt{\lambda}$, we have 
$$
P_1=aP+{b\over z}\Psi(P), \quad z_1=(a^2-b^2{\lambda})z=z.
$$
Since $a^2-({b\over z})^2Q(P)=a^2-b^2\lambda=1$, it follows from Corollary \ref{mufibres1} that
$\mu(P_1)=\mu(P)$. 

Conversely, if $\mu(P_1)=\mu(P)$ then, again by Corollary \ref{mufibres1}, there exist $x,y\in k$ such that $P_1=xP+y\Psi(P)$ and $x^2-y^2Q(P)=1$. If we set $u=x+yz\sqrt{\lambda}$, then
$u\bar{u}=1$  and $u\cdot(P,z)=(P_1,z_1)$.
\end{demo}

The above leads to natural action/angle variables in this context. For reasons of space and as they are not needed here, we shall omit them.

\subsection{Special symplectic $A_{\lambda}$-geometry}

Let $i:{\mathcal O}_{\lambda}\rightarrow V$ be inclusion.
We   can think of  $({\mathcal O}_{\lambda}, i)$ and $(\widehat{  {\mathcal O}_{\lambda} },\alpha)$ as respectively global and local geometric structures on ${\mathcal O}_{\lambda}$. Alternatively, we can think of  $(\widehat{  {\mathcal O}_{\lambda} }, i\circ p_{\lambda})$ and $(\widehat{  {\mathcal O}_{\lambda} },\alpha)$ as  global  geometric structures on $\widehat{  {\mathcal O}_{\lambda} }$ and then the `charts' $i\circ p_{\lambda}$ and $\alpha$ will each define distinguished classes of regular functions on $\widehat{  {\mathcal O}_{\lambda} }$ by pullback.

The first class consists of the regular  $k$-valued functions on $\widehat{  {\mathcal O}_{\lambda} }$  which  are the pullbacks by $ i\circ p_{\lambda}$  of linear functions on $V$.   If $k=\bb R$ this has a differential geometric interpretation : there is a unique torsion free, flat connection in the tangent bundle of $\widehat{  {\mathcal O}_{\lambda} }$ for which the exterior derivatives of these functions are covariantly constant.

The second class  is a  class of regular  $A_{\lambda}$-valued functions on $\widehat{  {\mathcal O}_{\lambda} }$ (i.e., elements of  $R(\widehat{  {\mathcal O}_{\lambda} })\otimes_kA_{\lambda}$). It consisits of  those functions which  are the pullbacks by $\alpha$ of  $A_{\lambda}$-linear functions on $V_{\lambda}$.  Again, if $k=\bb R$ this has a differential geometric interpretation but there are two cases.   If $\lambda$ is a square then $A_{\lambda}\cong \bb R\times\bb R$ and  $\widehat{  {\mathcal O}_{\lambda} }$ has a unique integrable paracomplex structure  for which these functions are paraholomorphic.
 If $\lambda$ is not a square  then $A_{\lambda}\cong \bb C$ and  $\widehat{  {\mathcal O}_{\lambda} }$ has a unique integrable complex structure  for which these functions are holomorphic.
 
  The two geometric structures above are compatible in the sense that the following diagram is commutative:
 \begin{equation}\label{}
\xymatrix
{
\widehat{  {\mathcal O}_{\lambda} }\ar[r]^{\alpha}\ar[d]_{p_{\lambda}}&V_{\lambda}\ar[d]^{{\mathcal R}e}\\
{\mathcal O}_{\lambda}\ar[r]^{i}&V.\\
}
\end{equation}
If $k=\bb R$ this also has a differential geometric interpretation: the connection $\nabla$ and (para)complex stucture $J$ satisfy $d^{\nabla}J=0$ where $d^{\nabla}$ is the exterior covariant derivative and $J$ is viewed as a 1-form with values in the tangent bundle. A real manifold with
a flat connection and (para)complex structure satisfying this equation is known as a special, (para)complex manifold [ACD].

\section{Appendix}
We describe here the SSRs which correspond  to simple Heisenberg graded  classical Lie algebras  by Theorem \ref{SSRandsimpleLiealg}.

\begin{ex} Let $(V,\omega)$ be a symplectic vector space of dimension $2n$.  We set ${\mathfrak m}=sp(V,\omega)$ and define $\tau:V\rightarrow sp(V,\omega)$ by
$$
{\tau}(v)(w)=\omega(v,w)v.
$$
The associated symmetric bilinear form $B_{\tau}$ tautologically satisfies   \eqref{ssrdefi}, and the Lie algebra ${\mathfrak g}({\mathfrak m},V,\omega,B_{\tau})$ is isomorphic to the symplectic Lie algebra $C_{n+1}$ of dimension $(n+1)(2n+3)$. The
cubic and quartic invariants vanish identically.
\end{ex}

\begin{ex}  \label{A_nSSR}
Let $(V,\omega)$ be a symplectic vector space of dimension $2n$ over $k$ and
let $J:V\rightarrow V$ be such that for some $\lambda\in k^*$,
\begin{align}
\nonumber
J^2&=\lambda Id,\\
\nonumber
\omega(J(v),w)+\omega(v,J(w))&=0\quad\forall v,w\in V.
\end{align}
Let ${\mathfrak m}$ be the commutant of  $J$ in $sp(V,\omega)$. The map 
$\mu:V\rightarrow sp(V,\omega)$ defined by
$$
\mu(v)=\tau(v)-{1\over \lambda}\tau(J(v))+{1\over 2\lambda}\omega(v,J(v))J
$$
takes its values in ${\mathfrak m}$ and the associated symmetric bilinear form $B_{\mu}$ satisfies  \eqref{ssrdefi},.  The cubic and normalised quartic covariants are
$$
\begin{cases}
\Psi(v)={3\over 2\lambda}\omega(v,J(v))J(v)\\
Q(v)={9\over 4\lambda}\omega(v,J(v))^2
\end{cases}
$$
If $\lambda$ is a square in $k$ then  the Lie algebra ${\mathfrak g}({\mathfrak m},V,\omega,B_{\mu})$ is isomorphic to  $A_{n+1}$. However if $\lambda$ is not a square in $k$, this is no longer true.
For example if $k=\bb R$ and $\lambda=-1$,  the isomorphism class of ${\mathfrak g}({\mathfrak m},V,\omega,B_{\mu})$ depends on the signature of the real quadratic form $q(v)=\omega(v,J(v))$: if $q$ is of signature $(2p,2(n-p))$, then ${\mathfrak g}({\mathfrak m},V,\omega,B_{\mu})$ is isomorphic to  $su(p+1,n-p+1)$ and ${\mathfrak m}$ is isomorphic to $u(p,n-p)$.

Returning to the general case, one can show that 
$\mu(v)=0$ iff $v$ is an eigenvector of $J$. Hence if $\lambda$ is not a square in
$k$, the only solution of $\mu(v)=0$ is $v=0$.
\end{ex}
\begin{ex} 
Let $(E,\Omega)$ be a two dimensional symplectic vector space and let $(F,g)$ be a nondegenerate $n$-dimensional quadratic space over $k$. Set
$$
V={\rm Hom}(E,F),\quad{\mathfrak m}={\mathfrak sl}(E,\Omega)\oplus \mathfrak {so}(F,g)
$$ 
and let ${\mathfrak m}$ act on $V$ by
$$
(s_1,s_2)\cdot A= s_2A-As_1
$$
If $A\in V$ we define its adjoint  $A^*\in {\rm Hom}(F,E)$  by
$$
\Omega(w,A^{*}(v))=g(A(w),v)\quad\forall v\in E,\forall w\in F.
$$
One checks that
 $(A^*B)^{t_{\Omega}}=-B^*A$, $(AB^*)^{t_{g}}=-BA^*$, $(AB^*C)^*=C^*BA^*$ and that the equation
$$
A^*B-B^*A=\Omega(A,B)Id_{E}
$$
defines an ${\mathfrak m}$-invariant symplectic form on $V$. 
If  $\mu:V\rightarrow{\mathfrak m}$ is defined by
$$
{\mu}(A)=\left(-A^*A, 2AA^*\right),
$$
the associated symmetric bilinear map satisfies   \eqref{ssrdefi} and the  Lie algebra ${\mathfrak g}({\mathfrak m},V,\omega,B_{\mu})$ is isomorphic to the orthogonal Lie algebra
$so(g\oplus2H)$ ($H$ is the hyperbolic plane). The cubic and normalized quartic covariants are
$$
\begin{cases}
\Psi(A)=3AA^*A,\\
Q(A)=-9{\rm det}(A^*A).
\end{cases}
$$
This example  becomes perhaps less obscure if we choose a basis  $\{e_1,e_2\}$  of $E$ such that $\Omega(e_1,e_2)=1$.  Then,  setting  $A(e_i)=a_i$ and $B(e_i)=b_i$ for $i=1,2$, we have
\begin{align}
\nonumber
\omega(A,B)&=g(a_1,b_2)-g(a_2,b_1),\\
\nonumber
\mu(A)(e,f)&=(-i_{e}(e_1\wedge e_2),2i_{f}(a_1\wedge a_2))\\
\nonumber
\Psi(A)(e)&=3i_{A(e)}(a_1\wedge a_2)\\
Q(A)&=9\left(g(a_1,a_2)^2-g(a_1,a_1)g(a_2,a_2)\right)
\end{align}
where $e\in E, f\in F$ and $i_e,i_f$  denote interior products with respect to $A^{\star}g$ and $g$ respectively.
\end{ex}

\end{document}